\documentclass[11pt]{amsart}
\usepackage{amsmath}
\usepackage{amsxtra}
\usepackage{amscd}
\usepackage{amsthm}
\usepackage{amsfonts}
\usepackage{amssymb}
\usepackage{eucal}

\newcommand{\nc}{\newcommand}
\newcommand{\rc}{\renewcommand}

\numberwithin{equation}{section}

\swapnumbers

\newtheorem{thm}{Theorem} [section]
\newtheorem{prop}[thm]{Proposition}
\newtheorem{lem}[thm]{Lemma}
\newtheorem{cor}[thm]{Corollary}
\newtheorem{rmk}[thm]{Remark}
\newtheorem{conj}[thm]{Conjecture}
\nc{\on}{\operatorname}

\nc{\bA}{{\mathbb A}}
\nc{\bB}{{\mathbb B}}
\nc{\bC}{{\mathbb C}}
\nc{\bD}{{\mathbb D}}
\nc{\bE}{{\mathbb E}}
\nc{\bF}{{\mathbb F}}
\nc{\bG}{{\mathbb G}}
\nc{\bH}{{\mathbb H}}
\nc{\bI}{{\mathbb I}}
\nc{\bJ}{{\mathbb J}}
\nc{\bK}{{\mathbb K}}
\nc{\bL}{{\mathbb L}}
\nc{\bM}{{\mathbb M}}
\nc{\bN}{{\mathbb N}}
\nc{\bO}{{\mathbb O}}
\nc{\bP}{{\mathbb P}}
\nc{\bQ}{{\mathbb Q}}
\nc{\bR}{{\mathbb R}}
\nc{\bS}{{\mathbb S}}
\nc{\bT}{{\mathbb T}}
\nc{\bU}{{\mathbb U}}
\nc{\bV}{{\mathbb V}}
\nc{\bW}{{\mathbb W}}
\nc{\bZ}{{\mathbb Z}}
\nc{\bX}{{\mathbb X}}
\nc{\bY}{{\mathbb Y}}

\nc{\cA}{{\mathcal A}}
\nc{\cB}{{\mathcal B}}
\nc{\cC}{{\mathcal C}}
\nc{\cD}{{\mathcal D}}
\nc{\cE}{{\mathcal E}}
\nc{\cF}{{\mathcal F}}
\nc{\cH}{{\mathcal H}}
\nc{\cI}{{\mathcal I}}
\nc{\cJ}{{\mathcal J}}
\nc{\cK}{{\mathcal K}}
\nc{\cL}{{\mathcal L}}
\nc{\cM}{{\mathcal M}}
\nc{\cN}{{\mathcal N}}
\nc{\cO}{{\mathcal O}}
\nc{\cP}{{\mathcal P}}
\nc{\cQ}{{\mathcal Q}}
\nc{\cR}{{\mathcal R}}
\nc{\cS}{{\mathcal S}}
\nc{\cT}{{\mathcal T}}
\nc{\cU}{{\mathcal U}}
\nc{\cV}{{\mathcal V}}
\nc{\cW}{{\mathcal W}}
\nc{\cZ}{{\mathcal Z}}
\nc{\cX}{{\mathcal X}}
\nc{\cY}{{\mathcal Y}}

\nc{\fA}{{\mathfrak A}}
\nc{\fB}{{\mathfrak B}}
\nc{\fC}{{\mathfrak C}}
\nc{\fD}{{\mathfrak D}}
\nc{\fE}{{\mathfrak E}}
\nc{\fF}{{\mathfrak F}}
\nc{\fG}{{\mathfrak G}}
\nc{\fH}{{\mathfrak H}}
\nc{\fI}{{\mathfrak I}}
\nc{\fJ}{{\mathfrak J}}
\nc{\fK}{{\mathfrak K}}
\nc{\fL}{{\mathfrak L}}
\nc{\fM}{{\mathfrak M}}
\nc{\fN}{{\mathfrak N}}
\nc{\fO}{{\mathfrak O}}
\nc{\fP}{{\mathfrak P}}
\nc{\fQ}{{\mathfrak Q}}
\nc{\fR}{{\mathfrak R}}
\nc{\fS}{{\mathfrak S}}
\nc{\fT}{{\mathfrak T}}
\nc{\fU}{{\mathfrak U}}
\nc{\fV}{{\mathfrak V}}
\nc{\fW}{{\mathfrak W}}
\nc{\fZ}{{\mathfrak Z}}
\nc{\fX}{{\mathfrak X}}
\nc{\fY}{{\mathfrak Y}}
\nc{\fa}{{\mathfrak a}}
\nc{\fb}{{\mathfrak b}}
\nc{\fc}{{\mathfrak c}}
\nc{\fd}{{\mathfrak d}}
\nc{\fe}{{\mathfrak e}}
\nc{\ff}{{\mathfrak f}}
\nc{\fg}{{\mathfrak g}}
\nc{\fh}{{\mathfrak h}}
\nc{\fiI}{{\mathfrak i}}  
\nc{\ffi}{{\mathfrak i}}  
\nc{\fj}{{\mathfrak j}}
\nc{\fk}{{\mathfrak k}}
\nc{\fl}{{\mathfrak{l}}}
\nc{\fm}{{\mathfrak m}}
\nc{\fn}{{\mathfrak n}}
\nc{\fo}{{\mathfrak o}}
\nc{\fp}{{\mathfrak p}}
\nc{\fq}{{\mathfrak q}}
\nc{\fr}{{\mathfrak r}}
\nc{\fs}{{\mathfrak s}}
\nc{\ft}{{\mathfrak t}}
\nc{\fu}{{\mathfrak u}}
\nc{\fv}{{\mathfrak v}}
\nc{\fw}{{\mathfrak w}}
\nc{\fz}{{\mathfrak z}}
\nc{\fx}{{\mathfrak x}}
\nc{\fy}{{\mathfrak y}}

\nc{\oA}{{\operatorname A}}
\nc{\oB}{{\operatorname B}}
\nc{\oC}{{\operatorname C}}
\nc{\oD}{{\operatorname D}}
\nc{\oE}{{\operatorname E}}
\nc{\oF}{{\operatorname F}}
\nc{\oG}{{\operatorname G}}
\nc{\oH}{{\operatorname H}}
\nc{\oh}{{\operatorname H}}
\nc{\oI}{{\operatorname I}}
\nc{\oJ}{{\operatorname J}}
\nc{\oK}{{\operatorname K}}
\nc{\oL}{{\operatorname L}}
\nc{\oM}{{\operatorname M}}
\nc{\oN}{{\operatorname N}}
\nc{\oO}{{\operatorname O}}
\nc{\op}{{\operatorname P}}
\nc{\oP}{{\operatorname P}}
\nc{\oQ}{{\operatorname Q}}
\nc{\oR}{{\operatorname R}}
\nc{\oS}{{\operatorname S}}
\nc{\oT}{{\operatorname T}}
\nc{\oU}{{\operatorname U}}
\nc{\oV}{{\operatorname V}}
\nc{\oW}{{\operatorname W}}
\nc{\oZ}{{\operatorname Z}}
\nc{\oX}{{\operatorname X}}
\nc{\oY}{{\operatorname Y}}

\rc{\bf}{\bfseries}
\nc{\noi}{{\noindent}}
\nc{\nop}{{\noindent {\bf Proof.}} }

\nc{\ra}{{	\rightarrow	}}

\nc{\df}{{ \aa{def}=		 }}		
\nc{\inv}{{ {}^{-1}      }}			

\nc{\lb}{\langle}             				
\nc{\rb}{\rangle}

\nc{\lB}{	\left(	}             			
\nc{\rB}{	\right)	}
\nc{\cBl}{{	\bbb{ \left( \right.}	}}             	
\nc{\cBr}{{	\bbb{ \left. \right)}	}}

\nc{\qlb}{{ \overline{{\mathbb Q}_l} }}      		
\nc{\ffq}{{  {\mathbb F}_q  }}           		
\nc{\ffp}{{  {\mathbb F}_p  }}           		
\nc{\tw}{   {}^{(1)}	}

\nc{\Spec}{{ \text{Spec}      		}}
\nc{\aand}{{\ \ \ \text{and}\ \ \ 	}}
\nc{\oor}{{\ \  \text{or}\ \ 	}}

\nc{\inj}{{	\hookrightarrow	}}    		
\nc{\injj}{{	\hookleftarrow	}}    		

\nc{\sur}{{	\twoheadrightarrow	}}	
\nc{\surr}{{	\twoheadleftarrow	}}	

\nc{\mm}{{\mapsto}}             		
\nc{\va}{{\uparrow}}              		

\nc{\G}{{	\mathcal G \text{r}		}}
\nc{\GG}{{	\G				}}
\def\GrG#1{{\mathcal G \text{r}^{#1}}}
\def\Gf#1{{\mathcal G \text{r}_{#1}}}

\nc{\GO} {{G_{{\mathcal O}}}}
\nc{\GK} {{{G_{{\mathcal K}}}}}
\nc{\go} {{\fg_{{\mathcal O}}}}
\nc{\gk} {{\fg_{{\mathcal K}}}}
\nc{\pg}{{	\operatorname P_{\GO}(\G,\Bbbk)		}}
\nc{\Hom}{{\operatorname{Hom}}}
\def\Gr#1{{\mathcal G \text{r}_{#1}}}

\nc{\cod}{{\operatorname{codim}}}
\nc{\ps}{{	\operatorname P_{\cS}	}}
\nc{\K}{{{\mathcal K}}}
\nc{\dime}{{\operatorname{dim}}}
\nc{\htt}{{\rho}}
\nc{\V}{{\operatorname{Mod}_\Bbbk}}
\def\pp#1{{{\operatorname P}_{\GO}(#1,\Bbbk)}}
\nc{\NK}{{N_\K}}
\nc{\NO}{{N_\cO}}
\nc{\ii}{{i\in I}}
\nc{\barr}{\overline}
\nc{\cdd}{{\cdot}}
\nc{\p}{{\on{P}}}
\nc{\dlm}{{{\tilde \G}^{\lambda,\mu}}}
\nc{\Glm}{{\GrG{\la+\mu}}}
\nc{\codim}{{\on{codim}}}
\nc{\Ext}{{\on{Ext}}}
\nc{\ext}{{\on{Ext}}}
\nc{\End}{\operatorname{End}}

\nc{\bb}{\underset}
\rc{\aa}{\overset}

\nc{\con}{{  \CD @>\cong>> \endCD 	}}  	
\nc{\conn}{{ \CD @<\cong<< \endCD  	}}  	

\nc{\ci}{{\circ}}               
\nc{\cd }{{\cdot}}              
\nc{\cddd}{{\cdot\cdot\cdot}}   

\nc{\tim}{{\times}}             
\nc{\ltim}{\ltimes}                  	%
\nc{\rtim}{\rtimes}			%

\nc{\ten}{{\otimes}}            
\nc{\bten}{{\boxtimes}}         
\nc{\pl}{{\oplus}}              

\nc{\al}{{\alpha }}
\nc{\be}{{\beta }}
\nc{\ga}{{\gamma }}
\nc{\de}{{\delta }}
\nc{\del}{{\partial }}
\nc{\ep}{{\varepsilon }}
\nc{\vap}{{\epsilon }}

\nc{\ze}{{\zeta }}
\nc{\et}{{\eta }}
\rc{\th}{{\theta }}

\nc{\vth}{{\vartheta }}

\nc{\io}{{\iota }}
\nc{\ka}{{\kappa }}
\nc{\la}{{\lambda }}
\nc{\vrho}{{\varrho}}
\nc{\si}{{\sigma }}
\nc{\ups}{{\upsilon }}
\nc{\vphi}{{\varphi }}
\nc{\om}{{\omega }}

\nc{\Ga}{{\Gamma }}
\nc{\De}{{\Delta }}
\nc{\nab}{{\nabla}}
\nc{\Th}{{\Theta }}
\nc{\La}{{\Lambda }}
\nc{\Si}{{\Sigma }}
\nc{\Ups}{{\Upsilon }}
\nc{\Om}{{\Omega }}

\nc{\A}{{\mathbb A }}
\nc{\B}{{\mathbb B}}
\nc{\C}{{\mathbb C}}
                \nc{\cc}{{\mathbb C}}
\nc{\Cs}{{\mathbb C^*}}
                \nc{\cs}{{\mathbb C^*}}
                \nc{\ccs}{{\mathbb C^*}}
\nc{\D}{{\mathbb D}}
\nc{\E}{{\mathbb E}}
\nc{\F}{{\mathbb F}}
        \nc{\hH}{{\mathbb H}}
\nc{\I}{{\mathbb I}}
\nc{\J}{{\mathbb J}}
        \nc{\lL}{{\mathbb L}}
\nc{\M}{{\mathbb M}}
\nc{\N}{{\mathbb N}}
        \nc{\pP}{{\mathbb P}}      
\nc{\Q}{{\mathbb Q}}
\nc{\R}{{\mathbb R}}
        \nc{\sS}{{\mathbb S}}
\nc{\T}{{\mathbb T}}
\nc{\U}{{\mathbb U}}
\nc{\W}{{\mathbb W}}
\nc{\Z}{{\mathbb Z}}
\nc{\X}{{\mathbb X}}
\nc{\Y}{{\mathbb Y}}

\rc{\k}{{\Bbbk}}   

\nc{\se}{       \section                }
\nc{\sus}{      \subsection             }
\nc{\sss}{      \subsubsection          }

\nc{\Lemm}{     \subsection{Lemma}              }
\nc{\lemm}{     \subsubsection{Lemma}           }
\nc{\slemm}{    \subsubsection*{Lemma}          }

\nc{\Pro}{      \subsection{Proposition}        }
\nc{\pro}{      \subsubsection{Proposition}     }
\nc{\spro}{     \subsubsection*{Proposition}    }

\nc{\Corr}{     \subsection{Corollary}          }
\nc{\corr}{     \subsubsection{Corollary}       }
\nc{\scorr}{    \subsubsection*{Corollary}      }

\nc{\Theo}{     \subsection{Theorem}            }
\nc{\theo}{     \subsubsection{Theorem}         }
\nc{\stheo}{    \subsubsection*{Theorem}                }

\nc{\remm}{     \subsubsection{Remarks}         }
\nc{\sremm}{    \subsubsection*{Remarks}                }
\nc{\rema}{     \subsubsection{Remarks}         }

\nc{\ex}{       \subsubsection{Example}         }
\nc{\sex}{      \subsubsection*{Example}                }
\nc{\exs}{      \subsubsection{Examples}                }
\nc{\sexs}{     \subsubsection*{Examples}               }

\nc{\lr}{{\leftrightarrow}}             
\nc{\lrs}{{\rightleftarrows}}           

\nc{\imp}{{\Rightarrow}}                
\nc{\eq}{{\Leftrightarrow}}             
        \nc{\Ra}{{\Rightarrow}}                 
        \nc{\LRa}{{\Leftrightarrow}}            

\nc{\nn}{{\newline}}
\nc{\indd}{{ ${} \ \ \ \ \  \ \ {} $	}}   
\nc{\inddd}{{ \indd\indd\indd		}}   

\nc{\bss}{{\backslash}}           

\nc{\ud}{       \underline      }               

\nc{\sub}{{     \subseteq       }}         
\nc{\subb}{{    \supseteq       }}         
\nc{\nsub}{{    \nsubseteq      }}         
\nc{\nsubb}{{   \nsupseteq      }}         %

\nc{\nin}{{     \notin  }}      

\nc{\ti}{\tilde}              
\nc{\tii}{\widetilde}         

\nc{\hatt}{\widehat}                            
\nc{\hata}{{    \mathbb{ \hat{} }          }}      

\nc{\ch}{\check}                                
\nc{\cha}{{     \mathbb{ \check{} }        }}      

\nc{\BBl}{{     \mathbb{ \left( \right.}   }}              
\nc{\BBr}{{     \mathbb{ \left. \right)}   }}

\nc{\h}{{       \hslash }}      

\nc{\All}{{     \forall }}
\nc{\Ex}{{      \exists         }}

\nc{\yy}{\infty}                        
\nc{\ys}{{  \frac{\infty}{2}  }}

\nc{\mat} {             \left(          \matrix }       
\nc{\emat}{             \endmatrix      \right) }
\nc{\sm} {              \left(          \smallmatrix    }       
\nc{\esm}{              \endsmallmatrix \right) }
\nc{\smat} {            \left(          \smallmatrix    }       
\nc{\esmat}{            \endsmallmatrix \right) }
\nc{\imat} {            \left.          \matrix }       
\nc{\eimat}{            \endmatrix      \right. }
\nc{\ism} {             \left.          \smallmatrix    }       
\nc{\eism}{             \endsmallmatrix \right. }

\nc{\ca}{               \left\{         \smallmatrix    }       
\nc{\eca}{              \endsmallmatrix \right\}        }
\nc{\Ca}{               \left\{         \matrix         }       
\nc{\Eca}{              \endmatrix      \right.         }

\rc{\AA}{{\mathcal A}}
\nc{\BB}{{\mathcal B}} 
\nc{\CC}{{\mathcal C}}
\nc{\DD}{{\mathcal D}}
\nc{\EE}{{\mathcal E}}
\nc{\FF}{{\mathcal F}}
\nc{\HH}{{\mathcal H}}
\nc{\II}{{\mathcal I}}
\nc{\JJ}{{\mathcal J}}
\nc{\KK}{{\mathcal K}}
\nc{\LL}{{\mathcal L}}
\nc{\MM}{{      \mathcal M        }}
\nc{\NN}{{\mathcal N}}
\nc{\OO}{{\mathcal O}}
\nc{\PP}{{\mathcal P}}
\nc{\QQ}{{\mathcal Q}}
\nc{\RR}{{\mathcal R}}
\rc{\SS}{{\mathcal S}}
\nc{\TT}{{\mathcal T}}
\nc{\UU}{{\mathcal U}}
\nc{\VV}{{\mathcal V}}
\nc{\WW}{{\mathcal W}}
\nc{\ZZ}{{\mathcal Z}}
\nc{\XX}{{\mathcal X}}
\nc{\YY}{{\mathcal Y}}

\nc{\px}{{	\PP_\SS(X,A)	}}
\nc{\py}{{	\PP_\TT(Y,A)	}}
\nc{\dx}{{	D_\SS(X,A)	}}
\nc{\dy}{{	D_\TT(Y,A)	}}

\nc{\Sn }{{  	S_\nu		}}
\nc{\Gl}{{  	\G^{\lambda}	}}
\nc{\Glb}{{  	\barr{\Gl}		}}

\nc{\xt}{{	X_*(T)		}}
\nc{\PGGk}{{ 	\PP_\GO(\GG) 	}}
\nc{\bast}{{	\framebox[2mm]{$\ast$}	}}

\nc\Ve{{	mod_\k^\ep		}}

\rc\cod{{\operatornamewithlimits{codim}}}
\nc{\Cal}{\cal}

\let\bc\C

\nc{\bbb}{ \boldsymbol         }  
\nc{\bul}{ \bullet         }  	
\nc{\rest}{{	\bbb{\big|}	}}
\nc{\tenZ}{{	\bb{\Z}\ten	}}
\nc{\tenk}{{	\bb{\k}\ten	}}

\nc{\tiii}{	\tii		}

\nc{\red}{ ^{red}	}
\nc{\thh}{ ^\text{th}	}
\nc{\af}{	^{aff}	}
\nc{\too}[1]{	\aa{#1}\to	}

\nc{\Gb}{  {\bbb G}	}

\nc{\Ker}{{	\text{Ker}	}}
\nc{\Coker}{{	\text{Coker}	}}
\rc{\Im}{{ 	\text{Im} 	}}
\nc{\rank}{{	\ \text{rank}\ 	}}
\nc{\Res}{{	\  \text{Res}   }}
\nc{\RHom}{{	\text{RHom}	}}
\nc{\HHom}{{	\text{$\HH$om}	}}
\nc{\EEnd}{{	\text{$\EE nd$}	}}
\nc{\AAut}{{	\text{$\AA ut$}	}}
\nc{\RHHom}{{	\text{R$\HH$om} }}

\let\cal\mathcal

\nc{\Gd}{{  {\check G} 		}}

\nc{\Gh}{{  \hat{\GG}		}}
\nc{\GA}{{  G(\AA)		}}
\nc{\BK}{{  B(\cal K)		}}
\nc{\BO}{{  B(\cal O)		}}
\nc{\BKz}{{  B(\cal K)_0	}}

\nc{\bk} {{\fb_{{\mathcal K}}}}

\nc{\tk}{{  	\frak b_\KK		}}
\nc{\nk}{{  	\frak n_\KK		}}

\nc{\Aut}{{\operatorname{Aut}}}

\nc{\limp}{{	\underset {\leftarrow}\lim	}}
\nc{\limi}{{	\underset {\rightarrow}\lim	}}

\nc{\supp}{{\operatorname{supp}}}

\nc{\PPP}{{	P	}}

\nc{\lab}{\label}

\nc{\phz}{{^p{\oh}^0}}
\nc{\Lbten}{{\aa{L}\boxtimes}}

\nc{\cG}{{	\check G	}}
\nc{\cGZ}{{	\check G_\Z	}}
\nc{\cGk}{{	\check G_\k	}}
\nc{\cGka}{{	\check G_\ka	}}
\nc{\cGQ}{{	\check G_\Q	}}

\nc{\cTZ}{{	\check T_\Z	}}
\nc{\cTk}{{	\check T_\k	}}
\nc{\cTka}{{	\check T_\ka	}}

\nc{\tG}{{	\tii G	}}
\nc{\tGZ}{{	\tii G_\Z	}}
\nc{\tGk}{{	\tii G_\k	}}
\nc{\tGka}{{	\tii G_\ka	}}
\nc{\tGQ}{{	\tii G_\Q	}}

\nc{\sG}{{	\tii G^*_\kappa	}}

\begin{document}

\title[Langlands duality and algebraic groups]{
Geometric Langlands duality and representations of algebraic groups over commutative rings}

\author{I. Mirkovi\'c}

\address {Department of mathematics, University of Massachusetts, Amherst, MA
01003, USA}
\email{ mirkovic@math.umass.edu}

\author{K. Vilonen}

\address{Department of Mathematics, Northwestern University, Evanston, IL
60208, USA}
\thanks{I. Mirkovi\'c and K. Vilonen were supported by NSF and the DARPA grant HR0011-04-1-0031}
\email {vilonen@math.northwestern.edu}

\maketitle

\nc{\uarr}{{\uparrow}}              		
\nc{\darr}{{\downarrow}}              		

\nc{\bn}{{	\fbox{\text{$\bbb{\ra}$}}	}}
\nc{\en}{{	\fbox{\text{$\bbb{<<}$}}	}}

\nc{\ftt}[1]{{\footnote{#1}}}
\nc{\fttt}[1]{{$^($\footnote{#1}$^)$}}
\nc{\bftt}[1]{\footnote{#1}}
\nc{\f}[1]{ \fbox{$ $}\footnote{ \fbox{!}#1 }\fbox{$ $}		}


\section{\bf Introduction}
\lab{Introduction}

In this paper we give a geometric version of the Satake isomorphism \cite{Sat}. As such,  it can be viewed as a first step in the geometric Langlands program. The connected complex reductive groups have a combinatorial classification by their root data. In the root datum the roots and the coroots appear in a symmetric manner and so the connected reductive algebraic groups come in pairs. If $G$ is a reductive group, we write $\check G$ for its companion and call it the dual group $G$. The notion of the dual group itself does not appear in Satake's paper,  but was introduced by Langlands, together with its various elaborations, in \cite{L1,L2} and is a corner stone of the Langlands program. It also appeared later in physics  \cite{MO,GNO}. In this paper we discuss the basic relationship between $G$ and $\check G$. 

We begin with a reductive $G$ and consider the affine Grassmannian  $\G$,  the Grassmannian for the loop group of $G$. For technical reason we work with formal algebraic loops. The affine Grassmannian is an infinite dimensional complex space. We consider a certain category of  sheaves, the spherical perverse sheaves, on $\G$. These sheaves can be multiplied using a convolution product and this leads to a rather explicit construction of a Hopf algebra, by what has come to be known as Tannakian formalism. The resulting Hopf algebra turns out to be the ring of functions on $\check G$. In this interpretation, the spherical perverse sheaves on the affine Grassmannian correspond to finite dimensional complex representations of $\check G$. Thus, instead of defining $\check G$ in terms of the classification of reductive groups, we provide a canonical construction of $\check G$, starting from $G$. We can carry out our construction over the integers. The spherical perverse sheaves are then those with integral coefficients, but the Grassmannian remains a complex algebraic object. The resulting $\check G$ turns out to be the Chevalley scheme over the integers, i.e., the unique split reductive group scheme whose root datum coincides with that of  the complex $\check G$. Thus, our result can also be viewed as providing  an explicit construction of the Chevalley scheme. Once we have a construction over the integers, we have one for every commutative ring and in particular for all fields. This provides another way of viewing our result: it provides a geometric interpretation of representation theory of algebraic groups over arbitrary rings. The change of rings on the representation theoretic side corresponds to change of coefficients of perverse sheaves, familiar from the universal coefficient theorem in algebraic topology. Note that for us it is crucial that we first prove our result for the integers (or $p$-adic integers) and then deduce the theorem for fields (of positive characteristic). We do not know how to argue the case of fields of positive characteristic directly. 

One of the key technical points of this paper is the construction of certain algebraic cycles that turn out to give a basis, even over the integers, of the cohomology of the standard sheaves on the affine Grassmannian. This result is new even over the complex numbers. These cycles are obtained by utilizing semi-infinite Schubert cells in  the affine Grassmannian. The semi-infinite Schubert cells can then be viewed as providing a perverse cell decomposition of the affine Grassmannian analogous to a cell decomposition for ordinary homology where the dimensions of all the cells have the same parity. The idea of searching for such a cell decomposition came from trying to find the analogues of the basic sets of \cite{GM} in our situation.

The first work in the direction of geometrizing the Satake isomorphism is \cite{Lu} where Lusztig  introduces the key notions and proves the result in the characteristic zero case on a combinatorial level  of affine Hecke algebras. Independently, Drinfeld had understood that geometrizing the Satake isomorphism is crucial for formulating the geometric Langlands correspondence. Following Drinfeld's suggestion, Ginzburg in \cite{Gi}, using \cite{Lu}, treated  the characteristic zero case of the geometric Satake isomorphism. Our paper is self-contained in that it does not rely on \cite{Lu} or \cite{Gi} and provides some improvements and precision even in the characteristic zero case. However, we make crucial use of an idea of Drinfeld, going back to around 1990. He discovered an elegant way of obtaining  the commutativity constraint by interpreting the convolution product of sheaves as a ``fusion" product. 

We now give a more precise version of our result. 
Let $G$ be a reductive algebraic group over the complex numbers. We write $\GO$ for the group scheme
$G(\bC[[z]])$ and $\G$ for the affine Grassmannian of $G(\bC((z)))/G(\bC[[z]])$; the affine Grassmannian is
an ind-scheme, i.e., a direct limit of schemes. Let $\k$ be a Noetherian, 
commutative unital ring of finite
global dimension. One can imagine $\k$ to be $\bC$, $\bZ$, or 
$\barr {\bF_q}$, for example. Let us write
$\pg$ for the category of
$\GO$-equivariant perverse sheaves with
$\k$-coefficients. Furthermore, let $\operatorname{Rep}_{\check G_\k}$ stand for the category
of 
$\k$-representations of
$\check G_\k$; here $\check G_\k$ denotes  the canonical smooth split reductive group scheme over $\k$
whose root datum is dual to that of $G$. The goal of this paper is to prove the following:
\begin{equation}
\text{the categories $\pg$ and $\operatorname{Rep}_{\check G_\k}$ are equivalent as tensor categories}\,.
\end{equation}
We do slightly more than this. We give a canonical construction 
of the group scheme $\check G_\k$ in
terms of $\pg$. In particular, we give a canonical construction of  the 
Chevalley group scheme $\check G_\bZ$ 
in terms of the complex group $G$.
This is one way to view our theorem. 
We can also view it as giving a geometric interpretation of representation
theory of algebraic groups over commutative rings. Although our results yield an interpretation of representation
theory over arbitrary commutative rings, note that on the geometric side we work over the complex numbers and
use the classical topology. The advantage of the  classical topology is that one can  work with sheaves with
coefficients in arbitrary commutative rings, in particular, we can use integer coefficients. Finally, our work can be
viewed as providing the unramified local geometric Langlands correspondence. 
In this context it is 
crucial that one
works on the geometric side also over fields other than $\bC$; this is easily done as the affine Grassmannian can be defined even   over the integers.
The modifications needed to do so are explained
in section \ref{Variants and the geometric Langlands program}. 
This can then be used to define
the notion of a Hecke eigensheaf
in the generality of
arbitrary systems of coefficients.

We describe the contents of the paper briefly. Section \ref{Perverse sheaves on the affine Grassmannian} is
devoted to the basic definitions involving the affine Grassmannian and the notion of perverse sheaves that we
adopt. In section
\ref{Semi-infiniteOrbits} we introduce our main tool, the weight functors. In this section we also give our crucial
dimension estimates, use them to prove the exactness of the weight functors, and, finally, we decompose the
global cohomology functor into a direct sum of the weight functors. The next section 
\ref{Convolution product} is devoted to putting a tensor structure on the category $\pg$; here, again, we make
use of the dimension estimates of the previous section. In section \ref{The commutativity constraint} we give,
using the Beilinson-Drinfeld Grassmannian, a commutativity constraint on the tensor structure. In section
\ref{Tensor functors} we show that global cohomology is a tensor functor and we also show that it is tensor
functor in the weighted sense. Section  \ref{The case of a field of characteristic zero} is devoted to the simpler
case when $\k$ is a field of characteristic zero. The next section \ref{Canonical basis} treats standard
sheaves and we show that their cohomology is given by specific algebraic 
cycles which provide a canonical basis for the cohomology. In the next section \ref{Representability of the weight functors} we prove that the weight
functors introduced in section \ref{Semi-infiniteOrbits} are representable. This, then,  will provide us with a
supply of projective objects. In section
\ref{The structure of  projectives that represent weight functors} we study the structure of these projectives
and prove that they have filtrations whose associated graded consists of standard sheaves. In section
\ref{Construction of the group scheme} we show that $\pg$ is equivalent, as a tensor category, to
$\operatorname{Rep}_{\tilde G_\k}$ for some group scheme $\tilde G_\k$. Then,  in the next section
\ref{identification}, we identify $\tilde G_\k$ with $\check G_\k$. A crucial ingredient in this section is the work
of Prasad and Yu \cite{PY}. We then briefly discuss in section
\ref{Representations of reductive groups} our results from the point of view of representation theory. In the final
section \ref{Variants and the geometric Langlands program} we briefly indicate how our arguments have to be
modified to work in the \'etale topology. 

Most of the results in this paper appeared in the announcement \cite{MiV2}. Since our announcement was published,  the papers \cite{Br} and \cite{Na} have  
appeared. Certain technical points that are necessary for us are treated 
in these papers. Instead of 
repeating the discussion here, we have chosen to refer to \cite{Br} and
\cite{Na} instead. Finally, let us note that we have not managed to carry out the idea of proof proposed in  \cite{MiV2} for theorem \ref{main result} (theorem 6.2 in \cite{MiV2}) and thus the paper \cite{MiV2} should be considered incomplete.   In this paper, as was mentioned above, we will appeal to \cite{PY} to prove theorem \ref{main result}.

We thank the MPI in Bonn, where some of this research was carried out. 
We also want to thank A .Beilinson, V. Drinfeld, and D. Nadler 
for many helpful discussions and KV wants to  thank G. Prasad and J. Yu
for answering a question  in the form of the paper \cite{PY}.

\section{\bf Perverse sheaves on the affine Grassmannian}
\lab{Perverse sheaves on the affine Grassmannian}

We begin this section by recalling the construction and the basic 
properties of the affine
Grassmannian $\G$. For proofs of these facts we refer to  \S4.5 of 
\cite{BD}. See also,
\cite{BL1} and \cite{BL2}.  Then we introduce the main object of 
study, the category
$\pg$ of equivariant perverse sheaves on $\G$.

Let $G$ be a complex, connected, reductive algebraic group.  We write ${\mathcal O}$ for
the  formal power
series ring $\bC[[z]]$ and ${\mathcal K}$ for its fraction field 
$\bC((z))$. Let
$G({\mathcal K})$ and $G({\mathcal O})$ denote, as usual, the sets of the
${\mathcal K}$-valued and the
${\mathcal O}$-valued points of $G$, respectively. The affine Grassmannian is
defined as  the quotient $G({\mathcal K})/G({\mathcal O})$.  The  sets
$G({\mathcal K})$ and $G({\mathcal O})$, and the quotient $G({\mathcal
K})/G({\mathcal O})$ have an algebraic structure over the complex numbers. The
space $G({\mathcal O})$ has a structure of a group scheme, denoted by
$\GO$, over $\bC$ and  the spaces
$G({\mathcal K})$ and $G({\mathcal K})/G({\mathcal
O})$ have structures of ind-schemes which we denote by $\GK$ and 
 $\G=\G_G$, respectively.
For us an ind-scheme means a direct limit of family of schemes where 
all the maps are
closed embeddings.  The morphism $\pi:\GK \to \G$ is locally trivial 
in the Zariski
topology, i.e., there exists a Zariski open subset $U\subset\G$ such that
$\pi^{-1}(U)\cong U\times\GO$ and $\pi$ restricted to $ U\times\GO$ is simply
projection to the first factor. For details see for example 
\cite{BL1,LS}. We write $\G$ as a
limit

\begin{equation} \label{limit}
\G \ = \ \varinjlim \Gf n\,,
\end{equation}
where the $\Gf n$ are finite dimensional schemes which are 
$\GO$-invariant. The group
$\GO$ acts on the $\Gf n$ via a finite dimensional quotient.

In this paper we consider
sheaves in the classical topology,
 with the exception of section \ref{Variants and the geometric Langlands program} where we use the etale topology. 
Therefore,  
it suffices for our purposes to consider the spaces
$\GO$,
$\GK$, and $\G$ as  reduced ind-schemes. We will do 
so for the rest of the
paper.

If $G=T$ is torus of rank $r$ then, as a reduced ind-scheme,
$\G \cong X_*(T) = \Hom (\bC^*,T)$, i.e., in this case the loop
Grassmannian is discrete. Note that, because $T$ is abelian, the loop
Grassmannian is a
group ind-scheme. Let $G$ be a reductive
group, write $Z(G)$ for the center of $G$ and let $Z = Z(G)^0$ denote
connected  component of the center. Let us further set 
$\barr G =
G/Z$. 
Then, as is easy to see, the map 
$\G_G\to \G_{\barr G}$
is a trivial covering with covering group $X_*(Z) 
= \Hom
(\bC^*,Z)$, i.e.,
$\G_G \cong \Gr{\barr G} \times 
X_*(Z)$, non-canonically. 
Note also that the connected components of $\G$ 
are exactly parameterized by the component group of $\GK$, i.e., by $\GK/(G_\K)^0$.
This latter group is isomorphic to $\pi_1(G)$, the topological
fundamental group of $G$.

The group scheme $\GO$ acts on $\G$ with finite dimensional orbits. In order
to describe the orbit structure, let us fix a maximal torus $T\subset G$. We
write $W$ for the Weyl group and $X_*(T)$  for the coweights
$\operatorname{Hom}(\bC^*,T)$. Then the $\GO$-orbits on $\G$ are parameterized
by the $W$-orbits in $X_*(T)$, and given $\lambda\in X_*(T)$ the $\GO$-orbit
associated to $W\la$ is 
$\GrG\lambda = \GO
\cdot L_\lambda\subset\G$, where $L_\lambda$ denotes the image of the 
point $\lambda
\in X_*(T)\ \sub\GK$ in $\G$. Note that the points $L_\lambda$ are precisely 
the $T$-fixed
points in the Grassmannian. To describe the closure relation between the $\GO$-orbits, we
choose a Borel $B \supset T$ and write $N$ for the 
unipotent radical of $B$. We use the convention that the roots in $B$ are the positive ones. 
Then, for dominant
$\lambda$ and $\mu$ we have
\begin{equation} \label{closure relation}
\GrG\mu\subset \overline{\GrG\lambda} \ \ \ \ \text{if and only if} 
\ \ \ \ \lambda-\mu \ \ \
\text{is a sum of positive coroots}\,.
\end{equation}

In a few arguments in this paper it will be important for us to 
consider a Kac-Moody
group associated to the loop group $\GK$.  Let us write 
$\De=\De(G,T)$ for 
the  root
system of $G$ with respect to $T$,  and we write similarly 
$\check\De=\check\De(G,T)$ for the coroots.
Let
$\Gamma \cong \bC^*$ denote the subgroup of automorphisms of $\K$ which acts
by multiplying the parameter
 $z\in\KK$ by $s\in   \bC^*\cong \Gamma$.
The group $\Gamma$ acts on $\GO$ and $\GK$ and hence we can form the
semi-direct product 
$\tii{G_\K} = \GK\rtimes\Gamma$. 
Then $\tii T = T
\times\Gamma$ is a Cartan subgroup  of $\tii {G_\K}$. 
An affine Kac-Moody group 
$\hatt{ G_\cK}$ is a
central extension, by the multiplicative group, of  $\tii{G_\K}$; note that 
the root systems are the same whether we consider $\tii {G_\K}$
or $\hatt{G_\cK}$. Let us write  $\delta\in X^*(\tii T)$ for 
the character   
which is trivial on $T$ and identity on
the factor $\Gamma \cong \bC^*$ and let
$\check\delta \in X_*(\tii T)$ 
be the
cocharacter
$\bC^*\cong\Gamma\subset
T\times \Gamma = \tii T$. We also view
the roots 
$\De$ as characters on 
$\tii T$, which are trivial on $\Gamma$. The $\tii T$-eigenspaces in $\gk$ are given by
\begin{equation}
(\gk)_{k \delta + \alpha}  \df
z^k \fg_\alpha  , \ k\in \bZ, \alpha \in 
\De\cup\{0\}\,,
\end{equation}
and thus the roots of $\GK$ are given by 
$\tii\De = \{\alpha+k\delta \in
X^*(\tii T)\mid \alpha\in \De\cup\{0\}, k\in \bZ\}-\{0\}$.

Furthermore, the orbit $G\cdot L_\lambda$ is isomorphic to the flag manifold
$G/P_\lambda$, where 
$P_\lambda$, the stabilizer of $L_\la$ in $G$, is a parabolic 
with a Levi factor associated to the roots
$\{\alpha\in\De\mid \lambda(\alpha)=0\}$. 
The orbit $\GrG\lambda$ can be viewed as a 
$G$-equivariant vector
bundle over $G/P_\lambda$. One way to see this is to observe that the 
varieties $G\cdot
L_\lambda$ are the fixed point sets of the $\bG_m$-action via the cocharacter
$\check\delta$. In this language,
\begin{equation}
\GrG\lambda \ = \ \{x\in\G \mid \lim_{s\to 0} \check\delta(s)x \in 
G\cdot L_\lambda\}
\end{equation}
In particular, the orbits
$\GrG\lambda$ are simply connected. If we choose a Borel $B$ 
containing $T$ and if we
choose the parameter
$\lambda \in X_*(T)$ of the orbit  $\GrG\lambda$ to be dominant, then
$\dim (\GrG\lambda) = 2\htt(\lambda)$, where 
$\rho\in X^*(T)$, as usual, is half the sum of positive roots with respect to $B$.
Let
us consider the  map
$\on{ev}_0:
\GO
\to G$, evaluation at zero. We write $I=\on{ev}_0\inv(B)$ for the 
Iwahori subgroup and
$K=\on{ev}_0\inv(1)$ 
for the highest congruence subgroup. The $I$-orbits are
parameterized by $X_*(T)$, and because the $I$-orbits are also 
$\on{ev}_0\inv(N)$-orbits, they are
affine spaces. This way each $\GO$-orbit acquires a cell decomposition as a union of $I$-orbits. The
$K$-orbit $K\cdot L_\lambda$ is the fiber of the vector bundle 
$\GrG\la \to G/P_\la$. Let
us consider the 
subgroup  ind-scheme
$G_\cO^-$ of $\GK$ whose
$\bC$-points consist of
$G(\bC[z \inv])$. The $G_\cO^-$-orbits are also indexed by $W$-orbits 
in $X_*(T)$
and the orbit attached to $\la\in X_*(T)$ is $G_\cO^-\cdot L_\la$. The
$G_\cO^-$-orbits are opposite to the
$\GO$-orbits in the following sense:
\begin{equation}
G_\cO^-\cdot L_\la \ = \ \{x\in\G \mid \lim_{s\to \infty} 
\check\delta(s)x \in G\cdot
L_\lambda\}\,.
\end{equation}
The above description implies that
\begin{equation}\label{opposite orbits}
(G_\cO^-\cdot L_\la)\cap \overline{\GrG\la} \ = \ G\cdot L_\la
\end{equation}
The group $G_\cO^-$ contains a negative level congruence subgroup
$K_-$ which is the kernel of the evaluation map $G(\bC[z \inv])\to G$ 
at infinity. Just as
for $\GO$, the fiber of the projection $G_\cO^-\cdot L_\la \to 
G/P_\la$ is $K_-\cdot
L_\la$.

We will recall briefly the notion of perverse sheaves that we will 
use in this paper
\cite{BBD}.
Let $X$ be a complex algebraic variety with a fixed (Whitney)
stratification
$\cS$.  
We also fix a commutative, unital ring $\k$. For simplicity of exposition we assume that $\k$ is
Noetherian of finite global dimension. This has the advantage of allowing us to work with finite complexes and finitely generated modules instead of having to use more complicated notions of finiteness. With suitable modifications, the results of this paper hold for arbitrary $\k$.   We denote by
$\operatorname{D}_{\cS}(X,\k)$ the bounded
$\cS$-constructible derived category of $\k$-sheaves. This is the full
subcategory of the derived category of $\k$-sheaves on $X$ whose objects
$\cF$ satisfy the following two conditions:
$$
\aligned i)\ \  & \oh^k(X,\cF) = 0 \ \ \ \text{for} \ |k|>0\,,
\\ ii)\ \ & \oh^k(\cF) \rest  S \ \ \ \  \text{is a local system of finitely
generated
$\k$-modules}
\\
&\text{ for all $S\in\cS$}\,.
\endaligned
$$
As usual we define the full subcategory $\ps(X,\k)$ of perverse
sheaves as follows. An $\cF\in\operatorname{D}_{\cS}(X,\k)$ is 
perverse if the following
two conditions are satisfied:
$$
\aligned i)\ \  & \oh^k(i^*\cF) = 0 \ \ \ \text{for} \ k> -
\dim_\bC S \
\text{ for any } \ i: S  \hookrightarrow X\,, S\in \cS\,,
\\ ii)\ \ & \oh^k(i^!\cF) = 0 \ \ \ \text{for} \ k< - \dim_\bC S
\
\text{ for any } \ i: S  \hookrightarrow X\,, S\in \cS\,.
\endaligned
$$
As is explained in \cite{BBD}, perverse sheaves $\ps(X,\k)$ form an
abelian category and there is a cohomological functor
$$
^p\oh^0 : \operatorname{D}_{\cS}(X,\k) \to \ps(X,\k)\,.
$$
Given a stratum $S\in \cS$ and $M$ a finitely generated $\k$-module then
$Rj_*M$ and $j_!M$ belong in $\operatorname{D}_{\cS}(X,\k)$. 
Following \cite{BBD} we
write
$^pj_*M$ for $^p\oh^0(Rj_*M) \in \ps(X,\k)$ and $^pj_!M$ for
$^p\oh^0(j_!M) \in\ps(X,\k)$. We use this type of notation systematically
throughout the paper. If $Y\subset X$ is locally closed and is a union
of strata in $\cS$ then, by abuse of notation, we denote by $\ps(Y,\k)$
the category $\operatorname{P}_{\mathcal T}(Y,\k)$, where $\mathcal T 
= \{S\in\cS \ \rest  \ S
\subset Y\}$.

Let us now assume that we have an action of a connected algebraic group $K$ on
$X$, given by $a:K\times X \to X$. Fix a Whitney
stratification $\cS$ of $X$ such that the action of $K$ preserves the
strata. Recall that an
$\cF\in\ps(X,\k)$ is said to be $K$-equivariant if there exists an
isomorphism
$\phi:a^*\cF \cong p^*\cF$ such that $\phi\rest \{1\}\times X = \text{id}$.
Here $p:K\times X \to X$ is the projection to the second factor. If
such an isomorphism $\phi$ exists it is unique. We denote by 
$\operatorname{P}_K
(X,\k)$ the full subcategory of $\ps(X,\k)$ consisting of equivariant
perverse sheaves. In a few instances we also make use of the equivariant derived category
$\operatorname{D}_{K}(X,\k)$, see \cite{BL}.

Let us now return to our situation. Denote the stratification 
induced by the
$\GO$-orbits on the Grassmannian $\G$  by $\cS$. The closed 
embeddings  $\Gf n \subset
\Gf m$, for $n\leq m$ induce embeddings of categories
$\operatorname {P}_{\GO}(\Gf n, \k) \to \operatorname 
{P}_{\GO}(\Gf m, \k)$.
This allows us to define the category of $\GO$-equivariant perverse 
sheaves on $\G$ as
$$
\pg \  =_{\text{def}} \ \varinjlim \operatorname{P}_{\GO}(\Gf n, \k)\,.
$$
Similarly we define $\ps (\G,\k)$, the category of perverse sheaves on
$\G$ which are constructible with respect to the $\GO$-orbits. In our
setting we have

\begin{prop}
\label{equivalence}
The categories $\ps (\G,\k)$ and 
$\pg$ are naturally
equivalent.
\end{prop}

We give a proof of this proposition in appendix
\ref{Categoryofequivariantperversesheaves}; the proof makes use of results of
section \ref{Semi-infiniteOrbits}.

Let us write $\operatorname{Aut}(\OO)$ for the group of automorphisms of the
formal disc $\operatorname{Spec}(\OO)$. The group scheme $\operatorname{Aut}(\OO)$ acts on
$\GK$, $\GO$, and $\G$. This action and the action of $\GO$ on the affine
Grassmannian extend to an action of the semidirect product $\GO \rtimes
\operatorname{Aut}(\OO)$ on $\G$. In the appendix
\ref{Categoryofequivariantperversesheaves} we also prove

\begin{prop}
\label{another equivalence}
The categories $\op_{\GO \rtimes \operatorname{Aut}(\OO)} (\G,\k)$ and 
$\pg$ are naturally
equivalent.
\end{prop}

\begin{rmk} If $\k$ is field of characteristic zero then
propositions \ref{equivalence} and
\ref{another equivalence} follow immediately from lemma \ref{semisimplicity}. 
\end{rmk}

Finally, we fix some notation that will be used throughout the paper. Given a $\GO$-orbit
$\GG_\la$, $\la\in X_*(T)$, and a $\k$-module $M$ we write
$\II_!(\la, M)$, $\II_*(\la,M)$, and $\II_{!*}(\la,M)$ for the perverse sheaves
${}^p\!j_!(M[\dim(\G^\la)])$, $j_{!*}(M[\dim(\G^\la)])$, and
${}^p\!j_*(M[\dim(\G^\la)])$, respectively; here $j:\G^\la \to \G$ denotes the inclusion.

\section{\bf Semi-infinite orbits and weight functors}
\label{Semi-infiniteOrbits}

Here we show that the global cohomology is a fiber functor for our 
tensor category. For $\k=\C$ this is proved by Ginzburg \cite{Gi}
and was treated earlier  in   \cite{Lu}, on the level of dimensions 
(the dimension of the 
intersection cohomology is the same as the dimension of the
corresponding representation).

Recall that we have fixed a maximal torus $T$, a Borel $B \supset T$ and denoted by $N$  the unipotent
radical of $B$. Furthermore, we write 
$\NK$  for the  group  ind-subscheme
of $\GK$ whose $\bc$-points are $N(\K)$. The
$\NK$-orbits on
$\G$ are parameterized by $X_*(T)$; to each $\nu\in X_*(T) = \Hom(\bC^*,T)$
we associate the $\NK$-orbit $S_\nu =
\NK\cdot L_\nu$. 
Note that these
orbits are neither of finite dimension nor of finite codimension.  We 
view them as
ind-varieties, in particular, their intersection with any
$\overline{\GrG\lambda}$ is an algebraic variety. The following 
proposition gives the
basic properties of these orbits. Recall that for $\mu,\la\in X_*(T)$ 
we say that
$\mu\leq\la$ if $\lambda-\mu$ is a sum of positive coroots.

\begin{prop} 
\label{one equation}
We have

(a) 
$
\overline{S_\nu}
=\ 
\cup_{\eta\leq \nu}\ S_\eta
$.

(b) Inside $\overline{S_\nu}$, the boundary of $S_\nu$ is given by
a hyperplane section under an embedding of $\G$ in projective space.

\end{prop}

\begin{proof}

Because translation by elements in $T_\K$ is an 
automorphism of  the 
Grassmannian, it
suffices to prove the claim on the identity component of the 
Grassmannian. Hence, we may
assume that $G$ is simply connected.
In that case $G$ is a product of simple factors and we may then furthermore assume that
$G$ is simple and simply connected.

For a positive  coroot 
$\check\alpha$,      there is $T$-stable $\bP^1$  passing through 
$L_{\nu-\check\alpha}$ such that
the remaining $\bA^1$ lies in $S_\nu$, constructed as follows.
First observe that the one parameter subgroup
$U_\psi$ for an affine root
$\psi=\al+k\de$ fixes $L_\nu$
if $z^{k-\lb\al,\nu\rb}\fg_\al$ fixes $L_0$, i.e., if
$k\ge \lb\al,\nu\rb$.
So, for any integer
$k<\lb\al,\nu\rb$,\  
$(\fg_\KK)_\psi$ 
does not fix $L_\nu$,
but 
$(\fg_\KK)_{-\psi}$ does. 
We conclude that for the $SL_2$-subgroup generated by the
one parameter subgroups
$U_{\pm\psi}$ the orbit
through $L_\nu$ is a $\bP^1$  and that $
U_\psi\cdot L_\nu
\cong
\bA^1  
$ lies in $S_\nu$ since $\al>0$.  
The point at infinity is then 
$L_{s_\psi\nu}$
for the reflection
$s_\psi$ in the affine root $\psi$. 
For 
$k=\lb\al,\nu\rb-1$ this yields 
$L_{\nu-\check\alpha}$ as the point at infinity.
Hence $S_{\nu-\check\alpha}\sub\overline{S_\nu}$ for any 
positive coroot
$\check\alpha$ and therefore $\cup_{\eta\leq \nu}\ 
S_\eta\subset\overline{S_\nu}$.

To prove the rest of the proposition we embed the ind-variety $\G$ in 
an ind-projective
space $\bP(V)$ via an ample line bundle $\cL$ on $\G$. For simplicity 
we choose $\LL$ to be the positive
generator of the Picard group of $\G$.  
The action of $\GK$ on
$\G$ only extends to a projective action on the line bundle $\cL$. To 
get an action on
$\cL$ we must  pass to the Kac-Moody group 
$\hatt{ G_\cK}$ associated to $\GK$,
which was discussed in the previous section. The highest weight $\Lambda_0$ of the resulting representation
$V=\on{H}^0(\G,\cL)$ is zero on  $T$ and one on the central $\bG_m$.
Thus, we get a  $\GK$-equivariant  embedding $\Psi:\G \hookrightarrow
\bP(V)$ which maps $L_0$ to the highest weight line  $V_{\Lambda_0}$. In
particular, the $T$-weight of the line $\Psi(L_0)=V_{\Lambda_0}$ is zero.

We need a formula for the $T$-weight of the line
$\Psi(L_\nu) = 
\nu\cdot\Psi(L_0) = \nu\cdot
V_{\Lambda_0}$. Now, $\nu\cdot V_{\Lambda_0} = V_{\tii\nu\cdot \Lambda_0}$,
where $\tii \nu$ is any lift of the element $\nu\in X_*(T)$ to
$\hatt{ T_\K}$, the restriction to 
$\tii {T_\K}$
of the  central extension $\hatt{ G_\cK}$ of 
$\tii {G_\cK}$ by $\bG_m$. 
For $t\in T$,
\begin{equation}
\begin{split}
(\tii\nu\cdot \Lambda_0)(t)
&= \Lambda_0({\tii\nu}^{-1}t\tii\nu) =
\Lambda_0({\tii\nu}^{-1}t\tii\nu t^{-1})\,,
\end{split}
\end{equation}
since $\Lambda_0(t) = 1$.
The commutator $x,y\mapsto xyx^{-1}y^{-1}$ on 
$\hatt{ T_\K}$ descends 
to a 
pairing of $T_\K\times T_\K$  to the central $\bG_m$. The restriction 
of this pairing to 
$X_*(T) \times T \to \bG_m$, can be viewed as a homomorphism $\iota:X_*(T) \to
X^*(T)$, or, equivalently, as a bilinear form $(\ , \, )_*$ on $X_*(T)$.  Since
$\Lambda_0$ is identity on the central $\bG_m$ and since 
${\tii\nu}^{-1}t\tii\nu
t^{-1}\in\bG_m$, we see that
\begin{equation}
(\tii\nu\cdot \Lambda_0)(t) \ = \ {\tii\nu}^{-1}t\tii\nu t^{-1} \ 
= \ (\iota\nu)(t)^{-1}\,,
\end{equation}
i.e., $\tii\nu\cdot \Lambda_0 = -\iota\nu$ on $T$. We will now describe the
morphism $\iota$. 

The description of the central extension of 
$\tii{\fg_\K}$, corresponding to
$\hatt{ G_\cK}$, makes use of an invariant bilinear form $(\ , \, )$ on $\fg$, 
see, for example,
\cite{PS}. 
From the basic formula for the coadjoint action of  $\hatt{ G_\cK}$
(see, 
for example, \cite{PS}),
it is clear that
the form $(\ ,
\, )_*$ above is the restriction of $(\ , \, )$ to $\ft = \bC\otimes 
X_*(T)$.
The form $(\ , \, )$ is 
characterized by the property
that the  
corresponding bilinear form 
$(\ , \, )^*$ on $\ft^*$
satisfies $(\theta,\theta)^* =2$ for the longest root 
$\theta$. 
Now, 
for a root $\alpha\in\De$ we find that
\begin{equation}
\iota\check\alpha \ = \ \frac 2 {(\alpha,\alpha)^*} \alpha \ = \ \frac
{(\theta,\theta)^*}{(\alpha,\alpha)^*} \alpha\in \{1,2,3\}\cdot \alpha
\end{equation}
We conclude that 
$\io(\Z\ch\De)\cap \Z_+\De_+=\io(\Z_+\ch\De_+)$, i.e.,
\begin{equation}\label{iota}
\nu < \eta  \ \ \ \text{is equivalent to}\  \ \ \iota\nu < \iota\eta 
\qquad \text{for}\ \nu,\eta\in
X_*(T)\,.
\end{equation}

Let us write
$V_{>-\iota \nu}\sub
V_{\ge -\iota \nu}$ 
for the sum of all 
the $T$-weight
spaces of $V$ whose $T$-weight is bigger than (or equal to)
$-\iota\nu$. 
Clearly the central
extension of $\NK$ acts by increasing the $T$ weights, i.e., its 
action preserves the
subspaces
$V_{>-\iota \nu}$ and $V_{\ge-\iota \nu}$. This, together with 
\eqref{iota},  implies that
$\cup_{\eta \le \nu}\ S_\eta = \Psi\inv(\bP(V_{\ge-\iota\nu}))$. In 
particular, $\cup_{\eta \le
\nu}\ S_\eta $ is closed. This, with $\cup_{\eta\leq \nu}\ 
S_\eta\subset\overline{S_\nu}$,
implies that $\overline{S_\nu}=\ \cup_{\eta\leq \nu}\ S_\eta$, 
proving part (a) of the
proposition.

To prove part (b), we first observe that  $\cup_{\eta < \nu}\
S_\eta = \Psi\inv(\bP(V_{>-\iota \nu}))$. The line
$\Psi(L_\nu)$ lies in $ V_{\ge-\iota \nu}$ but not in $V_{>-\iota \nu}$.
Let us 
choose a linear form
$f$ on $V$ which is non-zero on the line $\Psi(L_\nu)$ and which 
vanishes on all
$T$-eigenspaces  whose eigenvalue is different from $-\iota \nu$. Let 
us write $H$ for the
hyperplane $\{f=0\}\sub V$. 
By construction, for  $v\in \Psi(L_\nu)$, and
any $n$ 
in the central extension of $\NK$, $nv\in\ \bC^*\cdot v+\ V_{>-\iota \nu}$.
So $v\neq 0$
implies
$f(nv)\neq 0$, and we see that
$S_\nu\cap H = \emptyset$. Since $\cup_{\eta < \nu}S_\eta\subset H$,  we
conclude that $\overline{S_\nu}\cap H = \cup_{\eta < \nu}S_\eta$, as required.

\end{proof}

Let us also consider the unipotent radical
$N^-$ of the Borel $B^-$ opposite to $B$. 
The $N^-_\K$-orbits on
$\G$ are again parameterized by $X_*(T)$: to each $\nu\in X_*(T)$ we 
associate the
orbit $T_\nu = N^-_\K \cdot L_\nu$\,. The orbits $S_\nu$ and 
$T_\nu$ intersect the
orbits $\GrG\lambda$ as follows:

\begin{thm} 
\label{dimension estimates}
We have
$$
\aligned
\text{a)} \ \ &\text{The intersection} 
\ \ S_\nu\cap\GrG\lambda \ 
\text{is non-empty
precisely when} 
\ \ L_\nu\in\overline{\GrG\lambda} \ \ \text{and then}
\\
&\text{ $S_\nu\cap\Glb $  is of pure
dimension}
\ \ \rho(\nu+\lambda)\,,
\ \text{if $\lambda$ is chosen dominant}.
\\
\text{b)} \ \ &\text{The intersection} \ \ T_\nu\cap{\GrG\lambda} \ \ 
\text{is non-empty
precisely when} \ \ L_\nu\in\overline{\GrG\lambda} \ \ \text{and then}
\\
&\text{ $T_\nu\cap\barr{\GrG\lambda}$  is of pure
dimension}
\ \ -\rho(\nu+\lambda)\,,
\ \text{if $\lambda$ is chosen anti-dominant}\,.
\endaligned
$$
\end{thm}

\begin{rmk}
\label{weight picture}
Note that, by  \eqref{closure relation}, 
$L_\nu\in\overline{\GrG\lambda}$ if and only if
$\nu$ is a weight of 
the irreducible  representation of 
$\check G_\C$ of 
highest weight
$\lambda$; here $\check G$ is the complex Langlands dual group of $G$, i.e., the complex reductive group whose root datum is dual to that of $G$. 
\end{rmk}

\begin{proof}
It suffices to prove the statement a). Let  the
coweight $2\check\rho: \bG_m \to T$ be the sum of positive coroots. 
When we act by
conjugation by this coweight on
$\NK$, we see that for any element $n\in\NK$, $\lim_{s\to 0} 
2\check\rho(s) n = 1$.
Therefore any point
  $x\in S_\nu$ satisfies $\lim _{s\to 0} 2\check\rho(s) x = L_\nu$. As 
the $L_\nu$ are the
fixed points of the $\bG_m$-action via $2\check\rho$, we see that
\begin{equation}\label{S-orbit flow}
S_\nu \ = \ \{x\in\G \mid \lim _{s\to 0} 2\check\rho(s) x = L_\nu \}\,.
\end{equation}
Hence, if $x\in S_\nu\cap \GrG\lambda$ then, because $\GrG\lambda$ is 
$T$-invariant,
we see that $L_\nu\in\Glb$. Thus, $S_\nu\cap
\GrG\lambda$ is non-empty precisely when
$L_\nu\in\overline{\GrG\lambda}$. Recall that, as was remarked above, we then conclude, by \eqref{closure relation},  that  $S_\nu\cap
\GrG\lambda$ is non-empty precisely when $\nu$ is a weight of 
the irreducible  representation of 
$\check G_\C$ of 
highest weight
$\lambda$. 
Let us now assume that $\nu$ is such a wait.  

We begin with two extreme cases. We claim:
\begin{equation}
\label{extreme}
S_\nu\cap\overline{\GrG\nu} =\
N_\cO \cdot L_\nu= \begin{cases}
I\cdot L_\nu 	& \text{if $\nu$ is dominant}		\\
\{L_\nu\} 	& \text{if $\nu$ is anti-dominant}
\end{cases}
\end{equation}
We see this as follows. We first observe that $\NK = \NO \cdot 
(\NK\cap K_-)$. Then we
can write
\begin{equation}
S_\nu\cap\barr{\GrG\nu}=\
\NO\cdd (\NK\cap K_-)\cdd L_\nu\
\cap\ \barr{\GrG\nu}
\
=\
\NO\cdd \lB(\NK\cap K_-)\cdd\nu
\cap\barr{\GrG\nu}\rB\,.
\end{equation}
But now $(\NK\cap K_-)\cdd L_\nu\subset K_-\cdot L_\nu$ and by 
\eqref{opposite orbits}
we know that $G_\cO^-\cdot L_\nu\cap \overline{\GrG\lambda} = G\cdot L_\nu$ and
because  $K_- \cdot L_\nu$ is the fiber of the projection $G_\cO^-\cdot L_\nu
\to G\cdot L_\lambda$, we get $K_-\cdot \nu\cap 
\overline{\GrG\lambda} = L_\nu$. Thus
we have proved the first equality in \eqref{extreme}. If $\nu$ is 
antidominant, then $\NO$
stabilizes $L_\nu$. If $\nu$ is dominant then
$N_\cO^-$ stabilizes 
$L_\nu$ and then
$I\cdot L_\nu = B_\cO\cdot N_\cO^- \cdot L_\nu =B_\cO\cdot L_\nu = 
N_\cO\cdot L_\nu
$.

 From \eqref{extreme} we conclude that the theorem holds in the 
extreme cases when
$\nu=\lambda$ or $\nu = w_0\cdot \la$, where $w_0$ is the longest 
element in the Weyl
group. Let us now consider an arbitrary $\nu$ such that $L_\nu\in \overline{\GrG\la}$,
$\nu> w_0\cdot\la$ and let $C$ be an irreducible component of $S_\nu\cap \overline{\GrG\la}$. We will now relate this component to the two extremal cases above and make use proposition  \ref{one equation}.

Let us write $C_0$ for $C$, $d$ for the dimension of $C$,  and $H_\nu$ for the hyperplane of proposition  \ref{one equation} (b). Let us consider an irreducible component $D$ of $\bar C_0 \cap H_\nu$. By proposition \ref{one equation} the dimension of $D$ is $d -1$ and $D\subset\cup_{\mu<\nu} S_\mu$. Hence there is an $\nu_1<\nu=\nu_0$ such that $C_1 = D\cap S_{\nu_1}$ is open and dense in $D$. Of course $\dim C_1= d -1$. Continuing in this fashion we produce a sequence of coweights $\nu_k$, $k=0, \dots, d$, such that $\nu_k < \nu_{k-1}$, and a corresponding chain of irreducible components $C_k$ of $S_{\nu_k} \cap \overline{\GrG\la}$ such that $\dim C_k = d-k$. As dimension of $C_d$ is zero, we conclude that $\nu_d \geq w_0 \lambda$. Hence, we conclude that 
\begin{equation}
\label{dimension bound}
\dim C = d\leq \rho(\nu-w_0\cdot\lambda)\,.
\end{equation}

We now start from the opposite end. Let us write $A_0 = S_\lambda \cap  \overline{\GrG\la}$. Then, $\bar A_0 =  \overline{\GrG\la}$ and $\dim A_0 = 2\rho(\lambda)$. Let us proceed as before, and consider $\bar A_0\cap H_\lambda$. As $C\subset  \overline{\GrG\la}$, we can find a component $D$ of $ \bar A_0\cap H_\lambda$ such that $C\subset D$. Arguing just as above, there exists a $\mu<\la$ and a component $A_1$ of $S_\mu\cap  \overline{\GrG\la}$ such that $\bar A_1 = D$. Of course, $\dim A_1 = 2\rho(\lambda)-1$. Continuing in this manner  we can produce a  a sequence of coweights $\mu_k$, $k=0, \dots, e$, with $\mu_0=\lambda$, $\mu_e = \nu$, such that $\mu_k < \mu_{k-1}$, and a corresponding chain of irreducible components $A_k$ of $S_{\nu_k} \cap \overline{\GrG\la}$ such that $\dim A_k = 2\rho(\lambda)-k$ and $A_e=C$. From this we conclude that 
\begin{equation}
\label{codimension bound}
\codim_{\overline{\GrG\la}} C = e \leq \rho(\lambda-\nu)\,.
\end{equation}

The fact that 
\begin{equation}
\dim C + \codim_{\overline{\GrG\la}} C \ = \ \dim \overline{\GrG\la} \ = \ 2\rho(\la)\,,
\end{equation}
together with the estimates \eqref{dimension bound} and \eqref{codimension bound} force
\begin{equation}
\text{$\dim C = \rho(\nu-w_0\cdot\lambda)$ \ \ and \ \ $\codim_{\overline{\GrG\la}} C =  \rho(\lambda-\nu)$}\,,
\end{equation}
as was to be shown. 
\end{proof}

The corollary below will be used to construct the convolution 
operation on perverse sheaves
in the next section.

\begin{cor}\label{T-equivariant dimension estimate}
For any dominant $\lambda\in X_*(T)$ and any $T$-invariant closed 
subset $X\subset
\overline{\GrG\lambda}$ we have $\dim(X)\leq \max_{L_\nu\in
X^T}\ \htt(\lambda+\nu)$, where $X^T$ stands for the set of $T$--fixed 
points of $X$.
\end{cor}

\begin{proof}
 From the description \eqref{S-orbit flow} we see that $X\cap S_\nu$ 
is non-empty
precisely when $L_\nu\in X$. As
\begin{equation}
X \ = \ 
\cup_{L_\nu\in X^T}\ X\cap S_\nu\ \subset\ \cup_{L_\nu\in X^T}\
\overline{\G^\lambda} \cap S_\nu\,,
\end{equation}
we get our conclusion by appealing to  the previous theorem.
\end{proof}

Let us write $\V$ for the category of finitely generated $\k$-modules.

\begin{thm} 
\label{weightfunctors}
For all $\cA\in\pg$ we have a canonical isomorphism
\begin{equation}\label{updown}
\oh^k_c(S_\nu,\cA) \ 
\xrightarrow{\sim}\ 
\oh^k_{T_\nu}(\G,\cA)
\end{equation}
and both sides vanish for  $ k\neq 2\htt(\nu)\,.$

In particular, the functors
$
F_\nu:\ \pg\ra\ \V$,\ 
defined by 
$F_\nu\df\ \oh^{2\htt(\nu)}_c(S_\nu, -)
=$
$
\oh^{2\htt(\nu)}_{T_\nu}(\G,-)$,\
are exact.
\end{thm}

\begin{proof}
Let $\AA\in\pg$. For  any dominant $\eta$ the restriction 
$\AA\rest \GG^\eta$ lies, as a complex of sheaves, in degrees $\le
-\dim(\GG^\eta)=\ -2\htt(\eta)$, i.e., $\AA\rest \GG^\eta \in D^{\le
-2\htt(\eta)}(\GG^\eta,\k)$. From the dimension estimates
\ref{dimension estimates}  and fact that ${\oh}^{k}_{c}(S_\nu\cap\GG^\eta,\k)=0$, for 
$k>2\dim(S_\nu\cap\GG^\eta)=\ 
2\htt(\nu+\eta)$ we conclude: 
\begin{equation}
\label{NuEta}
\aligned
{\oh}^{k}_{c}(S_\nu\cap\GG^\eta,\cA)\  =\ 0 \ \ \ \text{if} \ k>2\htt(\nu)\,.
\endaligned
\end{equation}

A straightforward spectral sequence argument, filtering $\G$ by
$\overline{\G^\eta}$, implies that
${\oh}^{*}_{c}(S_\nu,\cA)$ can be expressed in terms of 
${\oh}^{*}_{c}(S_\nu\cap{\GrG\nu},\cA)$ and this implies  the first of the statements
below:
\begin{equation}\label{inequalities}
\aligned
&{\oh}^{k}_{c}(S_\nu,\cA)\  =\ 0 \ \ \ \text{if} \ k>2\htt(\nu)
\\
&{\oh}^k_{T_\nu}(\G,\cA) \ =\ 0 \ \ \ \text{if} \ k<2\htt(\nu)\,.
\endaligned
\end{equation}
The proof for  the second statement is completely analogous.

It remains to prove \eqref{updown}.  
Recall that we have a $\mathbb G_m$-action on $\G$ via the cocharacter 
$2\check\rho$ whose fixed points are the points $L_\nu$, $\nu\in X_*(T)$, and 
that 
\begin{align}
S_\nu \ &= \ \{x\in\G \mid \lim _{s\to 0} 2\check\rho(s) x = L_\nu \}
\\
T_\nu \ &= \ \{x\in\G \mid \lim _{s\to \infty} 2\check\rho(s) x = L_\nu \}\,.
\end{align}
The statement \eqref{updown} now follows from theorem 1 in \cite{Br}. 

\end{proof}

We will denote by
$F:\ \pg\ra\ \V$ 
the sum of the functors $F_\nu,\ \nu\in\ X_*(T)$.

\begin{thm} 
\label{exactness}
We have a natural equivalence of functors
$$
\mathbb {\oh}^* \ 
\cong \ F \ = \ \bigoplus_{\nu\in X_*(T)} \ 
\oh^{2\htt(\nu)}_c(S_\nu,-)\ :
\  \pg \to \V\,.
$$ 
Furthermore, the functors $F_\nu$ 
and this equivalence 
are independent of the choice of the
pair $T\subset B$. 
\end{thm}

\begin{proof} 
The Bruhat decomposition of
$G_\KK$
for the
Borel subgroups
$B_\KK,B^-_\KK$
gives decompositions
$\GG=\cup\ S_\nu=\cup \ T_\nu$ and hence two filtrations of $\GG$ by 
closures of $S_\nu$'s and $T_\nu$'s.
This gives  two filtrations of the
 cohomology functor $\mathbb {\oh}^*$, both  indexed by
$X_*(T)$. 
One is given by kernels of the 
morphisms of functors 
$\mathbb {\oh}^* \to
\oh^{*}_c( \barr {S_\nu},-)$ and the other 
by the images of
$\oh^{*}_{\barr {T_\nu}}(\G,-) \to \mathbb {\oh}^*$. 
The vanishing statement in
\ref{weightfunctors} implies that these filtrations are complementary. 
More precisely, in
degree $2\htt(\nu)$ we get $\oh^{2\htt(\nu)}_{\barr {T_\nu}}(\G,-) =
\oh^{2\htt(\nu)}_{T_\nu}(\G,-)$, $\oh^{2\htt(\nu)}_c(\barr {S_\nu},-) =
\oh^{2\htt(\nu)}_c(S_\nu,-)$, and the composition of the functors
$\oh^{2\htt(\nu)}_{T_\nu}(\G,-) \to \mathbb {\oh}^{2\htt(\nu)} \to
\oh^{2\htt(\nu)}_c(S_\nu,-)$ is the canonical equivalence in \ref{weightfunctors}. Hence, the
two filtrations of $\mathbb {\oh}^*$ split each other and provide the  desired natural
equivalence.

It remains to prove the independence of the equivalence and the functors $F_\nu$ of the
choice of $T\subset B$. Let us fix a reference $T_0\subset B_0$ and a $\nu\in X_*(T_0)$
which gives us the $S^0_\nu= (N_0)_\K \cdot\nu$. The choice of  pairs $T\subset B$ is
parameterized by the variety $G/T_0$. Note that there is a canonical isomorphism
between $T$ and $T_0$; they are both canonically isomorphic to the ``universal" Cartan
$B_0/N_0 = B/N$. Consider the following diagram
\begin{equation}
\begin{CD}
\G @<p<< \G \times G/T_0 @<j<< S
\\
@.  @V{q}VV @VrVV
\\
@. G/T_0 @= G/T_0\,.
\end{CD}
\end{equation}
Here $p,q,r$ are projections and $S=\{(x,gT_0)\in \G \times G/T_0 \mid x\in gS_\nu\}$. For a
point in $G/T_0$, i.e., for a choice of $T\subset B$ the fiber of $r$ is precisely the set
$S_\nu$ of the pair. Now, for any $\cA\in\pg$ the local system $Rq_*j_!j^*p^*\cA$ is a
sublocal system of $Rq_*p^*\cA$. As the latter local system is trivial, so is the former and
hence the functors $F_\nu$ are independent of the choice of $T\subset B$. 
\end{proof}

\begin{cor}\label{faithful and exact}
The global cohomology functor $\mathbb {\oh}^*=F: \pg
\to
\V$ is faithful and exact.
\end{cor}

\begin{proof}
The exactness follows from \ref{weightfunctors} and \ref{exactness}. If $\cA\in\pg$ is
non-zero then there exists an orbit $\G^\la$ which is open in the support of $\cA$. If we
choose $\la$ dominant then 
$T_\la\cap \barr{\G^\la}$
is a point in $\G^\la$ and we see that $F_\la(\cA)\neq
0$. As $\mathbb {\oh}^*$ does not annihilate non-zero objects it is faithful.
\end{proof}

\begin{rmk}
 The decompositions for $N$ and its opposite
unipotent subgroup  $N^-$  are explicitly
related by a canonical identification $\oh^k_{S_\nu}(\G,\cA) \cong
\oh^k_{T_{w_0\cdot\nu}}(\G,\cA)$, 
given by the action of
any representative of  $w_0$, the longest element in the Weyl group.
\end{rmk}

From the previous discussion we obtain the following criterion for a sheaf to be perverse:

\begin{lem}
\label{perversity criterion}
For a  sheaf $\cA\in\operatorname{D}_{\GO}(\G,\k)$, the following statements are
equivalent:
\begin{enumerate}
\item
The sheaf $\cA$ is perverse.
\item
For all $\nu\in X_*(T)$ the cohomology group
$\oh^{*}_c(S_\nu, \cA)$ is zero except possibly in degree $2\htt(\nu)$.
\item
For all $\nu$  group
$\oh^{*}_{S_\nu}(\G, \cA)$ is concentrated in degree $-2\htt(\nu)$.
\end{enumerate}
\end{lem}

\begin{proof}
By \ref{weightfunctors} and \ref{exactness} and an easy spectral sequence argument one
concludes that 
$\oh^{2\htt(\nu)}_c(S_\nu, \, ^p\!{\oh}^k(\cA)) = \oh^{2\htt(\nu)+k}_c(S_\nu, \,
\cA)$. This forces $\cA$ to be perverse. 
\end{proof}

Finally, we use the results of this section to give a rather explicit geometric
description of the cohomology of  the standard sheaves
$\II_!(\la,\k)$ and $ \II_*(\la,\k)$.

\begin{prop}
\lab{CanonicalBasis}
There  are 
canonical identifications
$$
F_\nu[\II_!(\la,\k)]\
\cong\
\k[Irr({\barr{\GG_\la}}\cap S_\nu)]
\
\cong\
F_\nu[\II_*(\la,\k)];
$$
here $\k[Irr({\barr{\GG_\la}}\cap S_\nu)]$ stands for the free $\k$-module generated by
the irreducible components of ${\barr{\GG_\la}}\cap S_\nu$.
\end{prop}

\begin{proof} We will give the argument for $\II_!(\la,\k)$. The argument for
$\II_*(\la,\k)$ is completely analogous. We proceed precisely the same way as in the
beginning of the  proof of \ref{weightfunctors}. Let us write
$\AA=\
\II_!(\la,\k)$. Consider 
 an orbit $\GG^\eta$  in the boundary of $\GG^\la$. Then
$\AA\rest \GG^\eta\in D^{\le -\dim(\GG^\eta)-2}(\GG^\eta,\k)$.
The estimate \ref{NuEta}  implies that ${\oh}^{k}_{c}(S_\nu\cap\GG^\eta,\cA)= 0$\
if   $k>2\htt(\nu)-2$. 
Therefore, we conclude by using the spectral sequence associated to the filtration of
$\G$ by $\overline{\G^\eta}$ that
${\oh}^{2\rho(\nu)}_{c}(S_\nu,\AA)
\cong\ {\oh}^{2\rho(\nu)}_{c}(S_\nu\cap\GG^\la,\AA)$.
Finally,
\begin{equation}
{\oh}^{2\rho(\nu)}_{c}(S_\nu\cap\GG^\la,\AA)
=\
{\oh}^{2\rho(\nu+\la)}_{c}(S_\nu\cap\GG^\la,\k)\ = \ 
{\oh}^{2\dim(S_\nu\cap\GG^\la)}_{c}(S_\nu\cap\GG^\la,\k)\,.
\end{equation}
As the last cohomology group is the top cohomology group, it is a free $\k$-module with
basis $Irr({\barr{\GG^\la}}\cap S_\nu)$.
\end{proof}

\section{\bf The
Convolution product}
\label{Convolution product}

In this section we will put a tensor category structure on $\pg$ 
via the convolution product. The idea that the convolution of perverse sheaves corresponds to 
tensor product of representations is due to Lusztig and the crucial proposition 4.2, 
for $\k = \bC$, is easy to extract from [Lu].
In some of our constructions in this section and the next one 
we are lead to sheaves with
infinite dimensional support. The fact that it is 
legitimate to work with such objects is explained in
section 2.2 of \cite{Na}. 

Consider the following diagram of maps
\begin{equation} \label{convolution}
\G\times\G \xleftarrow{p}\GK \times \G\xrightarrow{q} \GK \times_\GO
\G \xrightarrow{m} \G\,.
\end{equation}
Here $\GK \times_\GO\G$ denotes the quotient of $\GK \times \G$ by
$\GO$ where the action is given on the $\GK$-factor via right multiplication
by an inverse and on the $\G$-factor by left multiplication. The  $p$ and $q$
are projection maps and $m$ is the multiplication map.  We define the convolution
product
\begin{equation}\label{convolution product}
  \cA_1 \ast \cA_2 
\ = 
\ 
Rm_*\tii \cA \qquad\text{where \ $q^*\tii \cA  =
\ p^*(^p{\oh}^0(\cA_1\Lbten \cA_2))$}\,.
\end{equation}
To justify this definition, we note that the sheaf $p^*(^p{\oh}^0(\cA_1\Lbten \cA_2))$ on $\GK
\times\G$ is $\GO\times\GO$-equivariant with the first $\GO$ acting on the left and the second $\GO$
acting on the $\GK$-factor via right multiplication by an inverse and on the $\G$-factor by left
multiplication. As the second $\GO$-action is free, we see that the unique $\tii \cA$ in
\eqref{convolution} exists. 

\begin{lem}\label{outer tensor product is perverse}
If $\k$ is a field, or, more generally, if one of the factors ${\oh}^*(\G,\cA_i)$ is flat over $\k$, then the outer tensor product $\cA_1\Lbten \cA_2$ is perverse.
\end{lem}

When $\k$ is a field this is obvious on general grounds. When ${\oh}^*(\G,\cA_i)$ is flat over $\k$
one sees this by applying Lemma \ref{perversity criterion} to the Grassmannian $\GG\tim\GG$ 
of  $G\tim G$. First, as ${\oh}^*(\G,\cA_i)$ is flat, so are its direct summands $\oh^{2\htt(\nu_i)}_c(S_{\nu_i}, \cA_i)$. Now we have
\begin{equation}
\oh^{k}_c(S_{\nu_1}\times S_{\nu_2}, \cA_1\Lbten \cA_2) \ = \ \bigoplus_{k_1 + k_2 = k} \oh^{k_1}_c(S_{\nu_1}, \cA_1) \overset{L}\ten \oh^{k_2}_c(S_{\nu_2},  \cA_2)\,.
\end{equation}
By the flatness assumption the tensor product on the right has no derived functors. Hence, $\oh^{k}_c(S_{\nu_1}\times S_{\nu_2}, \cA_1\Lbten \cA_2) = 0$ if $k\neq 2\htt(\nu_1+\nu_2)$. Therefore, by Lemma \ref{perversity criterion}, $\cA_1\Lbten \cA_2$ is perverse.

\begin{prop}\label{convolution preserves perversity}
 The convolution product $  \cA_1 \ast \cA_2 $
of two perverse sheaves is perverse.
 \end{prop}

To prove this, let us introduce the notion of a stratified semi-small 
map. To this end,
let us consider two complex  stratified spaces $(Y,\cT)$ and $(X,\cS)$ and a
map $f:Y\to X$. We assume that the two stratifications are locally trivial with
connected strata and that $f$ is a stratified with respect to the
stratifications $\cT$ and $\cS$, i.e., that for any $T\in\cT$ the image
$f(T)$ is a union of strata in $\cS$ and for any $S\in\cS$ the map
$f\rest f^{-1}(S):  f^{-1}(S)
\to S$ is locally trivial in the stratified sense. We say that
$f$ is a stratified semi-small map if
\begin{equation}
\begin{aligned}
a) \ \ 		&	\text{for any $T\in\cT$ the map $f\rest \barr T$ is proper}
\\
b) \ \ 		&	\text{for any $T\in\cT$ and  any $S\in\cS$
				such that $S\subset	f(\barr T)$ we have}
\\
		&	\dim(f^{-1}(x)\cap \barr T) \leq \frac 1 2 
					(\dim f(\barr T) - \dim S)
\\
		&	\text{for any (and thus all) $x\in S$\,. }
\end{aligned}
\end{equation}
Let us also introduce the notion of a small stratified map. We say 
that $f$ is a small
stratified map if there exists a (non-trivial) open dense stratified subset $W$ of
$Y$ such that
\begin{equation}\label{small}
\begin{aligned}
  a) \ \ &\text{for any $T\in\cT$ the map $f\rest \barr T$ is proper}
\\ b) \ \ &\text{the map $f\rest W:W\to f(W)$ is finite and $W=f^{-1}(f(W))$}
\\ c) \ \ &\text{for any $T\in\cT$ and  any $S\in\cS$ such that
$S\subset f(\barr T)-f(W)$}
\\ &\text{we have}\ \ \dim(f^{-1}(x)\cap \barr T) < \frac 1 2 (\dim
f(\barr T)
-
\dim S)
\\ &\text{for any (and thus all) $x\in S$\,. }
\end{aligned}
\end{equation}

The result below follows directly from dimension counting:

\begin{lem}\label{small and semi-small maps}
If $f$ is a semi-small stratified map then
$Rf_*\cA\in\p_{\cS}(X,\k)$ for all $\cA\in\p_\cT(Y,\k)$\,. If $f$ is a small
stratified map then, with any $W$ as above, and any $\cA\in\p_\cT(W,\k)$, we
have $Rf_* j_{!*} \cA = \tii j_{!*} f_*\cA$, where $j: W \hookrightarrow Y$ and
$\tii j :f(W) \hookrightarrow X$ denote the two inclusions.
\end{lem}

We apply the above considerations, in the semi-small case,  to our situation.
We take $Y=\GK\times_\GO\G$ and  choose $\cT$ to be the stratification whose
strata are $\dlm=p^{-1}(\GrG\lambda)\times_\GO\GrG\mu$, for
$\lambda,\mu\in X_*(T)$\,. We also let $X=\G$, $\cS$ the stratification by
$\GO$-orbits, and choose
$f=m$. Note that the sheaf $\tii \cA$ is constructible with respect 
to the stratification
$\cT$.
To be able to apply \ref{small and semi-small maps} and 
conclude the  proof of
\ref{convolution preserves perversity},  we  appeal to the following

\begin{lem}
\lab{stratified semi-small}
The multiplication map
$\GK \times_\GO \G \xrightarrow {m} \G$ is a stratified semi-small 
map with respect to
the stratifications above.
\end{lem}

\begin{proof}
We need to check that for any
$\GO$-orbit $\GrG\nu$ in $\barr{\GrG{\la+\mu}}$,
the dimension of the fiber
$m\inv L_\nu\cap\dlm$ of $m:\dlm\ra\barr{\Glm}$ at $L_\nu$,
is not more than
$
{\frac{1}{2}}
\codim_{\barr\Glm}\G_\nu$.
We can  assume that
$\nu$ is anti-dominant since $\GrG{w\cdd \eta}=\GrG{\eta},\ w\in W.$
Since for any dominant $\eta$,
$\dim\GrG\eta=2\htt(\eta)$, the codimension in question is:
$$
\codim_{\barr\Glm}\GrG\nu=2\htt(\la+\mu)-2\htt(w_0\cdd \nu)=
2\htt(\la+\mu+\nu).
$$
Therefore, we need to show that
\begin{equation} \label{estimate}
\dim(m\inv L_\nu\cap\dlm)\le \
\htt(\la+\mu+\nu)\,.
\end{equation}

Let $p$ be the projection
$\GK\times_\GO \G \to \G$ given by $(g,h\GO)\mapsto g\GO$,
and consider the isomorphism
$(p,m):\ \GK\times_\GO \G \ \cong \ \G\times\G$.
The mapping $(p,m)$ carries the fiber $m\inv L_\nu$
to $\G\times L_\nu$. The set
$p(m\inv L_\nu\cap \barr\dlm)$ 
is $T$-invariant, and hence we can
apply corollary \ref{T-equivariant dimension estimate} to compute its 
dimension. To do
so, we need to find the $T$-fixed points in $p( m\inv L_\nu\cap\barr \dlm )\
\subset\
\barr{\GrG\la}$. The $T$-fixed points in $m\inv L_\nu\cap {\dlm}$ are 
precisely the
points $(z^\phi,z^\psi\GO)$ such that
$\phi$ and $\psi$ are weights of $L(\la)$ and $L(\mu)$ and $\phi+\psi=\
\nu$. Hence, the set $T$-fixed points in $m\inv L_\nu\cap 
\barr{\dlm}$ consists of the
points of the form
$(z^\phi,z^\psi\GO)$ with $\phi+\psi=\ \nu$
and
$\phi$ and $\psi$  weights of irreducible representations
$L(\la')$ and $L(\mu')$ for some
  dominant $\la',\mu'$ such that $\la'\le\la,\ \mu'\le\mu$.
For $\phi,\psi,\mu'$ as above, we have
$$
\htt(\la+\phi)\
\le\
\htt(\la+\phi) +\htt(\psi+\mu')
=\
\htt(\la+\nu+\mu')\le\
\htt(\la+\nu+\mu)\,.
$$
Therefore,
\begin{equation}
\dim(p(m\inv L_\nu\cap \barr\dlm)) \ \leq \ \max_{L_\phi\in
p(m\inv L_\nu\cap \barr\dlm)^T}(\htt(\lambda+\phi)) \leq \htt(\la+\nu+\mu)\,.
\end{equation}
This implies \eqref{estimate} and concludes the proof.

\end{proof}

\begin{rmk} One can also prove proposition
\eqref {convolution preserves perversity} and lemma 
\eqref{stratified semi-small} in a less geometric way by a rather direct translation of Lusztig's results on affine Hecke algebras in \cite{Lu}.
In the case of fields of characteristic zero, proposition
\eqref {convolution preserves perversity} was first proved by Ginzburg in \cite{Gi} in this manner.
From the characteristic zero case one can deduce
lemma \eqref{stratified semi-small}
and therefore also the general case of
\eqref {convolution preserves perversity}.
\end{rmk}

In completely analogy with \eqref{convolution product}, we can define 
directly the
convolution product of three sheaves, i.e., to $\cA_1,\cA_2,\cA_3$ we can associate a
perverse sheaf 
$\cA_1\ast\cA_2\ast\cA_3$. Furthermore, we get canonical isomorphisms
$\cA_1\ast\cA_2\ast\cA_3\cong (\cA_1\ast\cA_2)\ast\cA_3$ and $\cA_1\ast\cA_2\ast\cA_3\cong
\cA_1\ast(\cA_2\ast\cA_3)$. This yields a functorial isomorphism
$(\cA_1\ast\cA_2)\ast\cA_3\cong \cA_1\ast(\cA_2\ast\cA_3)$.

Thus we obtain:

\begin{prop} The convolution product \eqref{convolution product} on the abelian category $\pg$ is associative. 
\end{prop}

\section{\bf
The commutativity constraint and the fusion product}
\lab{The commutativity constraint}
In this section we show that the convolution product defined in the last section can be
viewed as a ``fusion" product. This interpretation allows one to provide the
convolution product on
$\pg$ with a commutativity constraint, making $\pg$ into an associative,
commutative tensor category.  The exposition follows very
closely that in \cite{MiV2}. The idea of interpreting the convolution product as a fusion
product and obtaining the commutativity constraint in this fashion is due to Drinfeld and was communicated to us by Beilinson.

Let $X$ be a smooth complex algebraic curve. For a closed point $x\in X$ we write
$\cO_x$  for the completion of the local ring at $x$ and  $\cK_x$ for its
fraction field. Furthermore, for a $\bC$-algebra $R$ we write $X_R =
X\times\text{Spec}(R)$, and $X^*_R= (X-\{x\})\times \text{Spec}(R)$\,. Using the results of
~\cite{BL1,BL2,LS} we can now view the Grassmannian
$\GG_ {x} = G_{\cK_x}/G_{\cO_x}$ in the following manner.  It is the ind-scheme which
represents the
 functor from
$\bC$-algebras to sets :
$$ R \mapsto \{ \cF \text{  a  $G$-torsor on $X_R$,
$\nu: G\times X^*_R \ra \cF\rest_{ X^*_R}$ a trivialization on $X^*_R$  }
\}\,.
$$
Here the pairs $(\cF,\nu)$ are to be taken up to isomorphism. 

Following \cite{BD} we  globalize this
construction and at the same time work over several copies of the curve.
Denote the
$n$ fold product by
$X^n = X\times
\dots\times X$  and consider the functor
\begin{equation}\label{global affine grassmannian}
R \mapsto \left\{\aligned &(x_1,\dots, x_n) \in X^n(R), \ \ \cF\text{ a
$G$-torsor on
$X_R$\,, }\\ &\text{$\nu_{(x_1,\dots,x_n)}$ a trivialization of
$\cF$ on $X_R - \cup {x_i}$}\endaligned\right\}\,.
\end{equation}
Here we think of the points $x_i: \text{Spec}(R) \to X$ as subschemes of
$X_R$ by taking their graphs. This functor  is
represented by an ind-scheme $\GG_ {X^{n}}$. Of course  $\GG_
  {X^{n}}$ is an
ind-scheme over
$X^n$ and its fiber over the point
$(x_1,\dots,x_n)$ is simply $\prod_{i=1}^k \GG_ {y_i}$\,, where
$\{y_1,\dots,y_k\}=\{x_1,\dots,x_n\}$, with all the $y_i$ distinct. We write
$\GG_ {X^{1} }= \GG_ X$.

We will now extend the diagram of maps \eqref{convolution}, which was 
used to define the
convolution product, to the global situation, i.e., to a 
diagram of ind-schemes
over $X^2$:
\begin{equation} \label{global convolution}
\GG_ X\times\GG_ X \xleftarrow{p} \widetilde{\GG_ X\times \GG_ X} \xrightarrow{q}
\GG_ X
\tii
\times
\GG_ X  
\xrightarrow{m} 
\GG_ {X^{2}}
\xrightarrow{\pi} 
{X^{2}}
\,.
\end{equation}

Here,
$\widetilde{\GG_ X\times \GG_ X}$ denotes the ind-scheme representing the functor
\begin{equation}\label{xxx}
R \mapsto \left\{\aligned &(x_1,x_2)\in X^2(R); \ \cF_1,\cF_2\text{
$G$-torsors on
$X_R$;\
$\nu_i$ a trivialization of }
\\ &\text{$\cF_i$ on $X_R -x_i$, for $i=1,2$; \ $\mu_1$ a trivialization of
$\cF_1$ on
$\widehat{(X_R)}_{x_2}$}
\endaligned
\right\},
\end{equation}
where $\widehat{(X_R)}_{x_2}$ denotes the formal neighborhood of
$x_2$ in
$X_R$. The \lq\lq twisted product" $\GG_ X  \tii \times \GG_ X $ is the
ind-scheme representing the functor
\begin{equation}
R \mapsto\left\{\aligned &(x_1,x_2)\in X^2(R); \ \cF_1,\cF \text{
$G$-torsors on $X_R$; $\nu_1$ a trivialization }
\\ &\text{of
$\cF_1$ on $X_R-x_1$;} \ \eta : \cF_1 \rest_{(X_R -x_2)} \xrightarrow{\ \simeq\ }
\cF\rest_{(X_R -x_2)}
\endaligned
\right\}\,.
\end{equation}
It remains to describe the morphisms $p$, $q$, and $m$ in 
\eqref{global convolution}.
Because all the spaces in \eqref{global convolution} are ind-schemes 
over $X^2$,  and all
the functors  involve the choice of the same point $(x_1,x_2)\in 
X^2(R)$, we omit it in the
formulas below. The morphism $p$ simply forgets the choice of 
$\mu_1$, the morphism
$q$ is given by the natural transformation
$$
  (\cF_1,\nu_1,\mu_1;\cF_2,\nu_2) \mapsto (\cF_1,\nu_1,\cF,\eta),
$$
where $\cF$ is the $G$-torsor gotten by gluing $\cF_1$ on $X_R - x_2$ and
$\cF_2$ on $\widehat{(X_R)}_{x_2}$ using the isomorphism induced by
$\nu_2\circ\mu_1^{-1}$ between $\cF_1$ and $\cF_2$ on
$(X_R-x_2)\cap \widehat{(X_R)}_{x_2}$. The morphism
$m$ is given by the natural transformation
$$
(\cF_1,\nu_1,\cF,\eta) \mapsto (\cF,\nu)\,,
$$
where $\nu = (\eta \circ \nu_1)\rest (X_R - x_1-x_2)$.

The global analogue of $\GO$ is the group-scheme 
$G_{X^{n},\mathcal O}$ which
represents the functor
\begin{equation}
R \mapsto \left\{\aligned &(x_1,\dots,x_n) \in X^n(R), \ \ \cF\text{ the
trivial
$G$-torsor on
$X_R$\,, }\\ &\text{$\mu_{(x_1,\dots,x_n)}$ a trivialization of
$\cF$ on $\widehat{(X_R)}_{(x_1\cup\dots\cup x_n)}$}\endaligned\right\}\,.
\end{equation}

Proceeding as in section \ref{Convolution product} we define the convolution
product of 
$\cB_1,\cB_2\in
\operatorname{P}_{G_{X, \mathcal O}}(\GG_ X,\k)$ 
by the formula
\begin{equation}\label{Global convolution product}
\cB_1 \ast_X \cB_2
\ = \ 
Rm_*\ 
\tii \cB 
\qquad
\text{where \ $q^*
\tii \cB = 
p^*(\phz(\cB_1\Lbten \cB_2))$}\,.
\end{equation}
To make sense of this definition, we have to explain how the group scheme ${G_{X, \mathcal O}}$ acts on various spaces. First, to see that it acts on $\GG_X$, we observe that  we can rewrite the functor in \eqref{global affine grassmannian}, when $n=1$, as follows:

\begin{equation}\label{alternate global affine grassmannian}
R \mapsto \left\{\aligned &x \in X(R), \ \ \cF\text{ a
$G$-torsor on
$\widehat{(X_R)}_x $\,, }\\ &\text{$\nu_x$ a trivialization of
$\cF$ on $\widehat{(X_R)}_x - {x}$}\endaligned\right\}\,.
\end{equation}

Thus we see that ${G_{X, \mathcal O}}$ acts on $\GG_X$ by altering the trivialization in \eqref{alternate global affine grassmannian} and hence we can define the category $\operatorname{P}_{G_{X, \mathcal O}}(\GG_ X,\k)$. As to $\widetilde{\GG_ X\times \GG_ X}$, two actions of  ${G_{X, \mathcal O}}$ are relevant to us. First, let us view ${G_{X, \mathcal O}}$ as group scheme on $X^2$ by pulling it back for the second factor. Then ${G_{X, \mathcal O}}$ acts by altering the trivialization $\mu_1$ in \eqref{xxx}. This action is free and exhibits $p: \widetilde{\GG_ X\times \GG_ X} \to \GG_X\times\GG_X$ as a ${G_{X, \mathcal O}}$ torsor. To describe the second action we rewrite the definition of  $\widetilde{\GG_ X\times \GG_ X}$ in the same fashion as we did for $\GG_X$, i.e., $\widetilde{\GG_ X\times \GG_ X}$  can also be viewed as representing the functor

\begin{equation}
R \mapsto \left\{\aligned &(x_1,x_2)\in X^2(R); \ \text{for $i=1,2$} \ \cF_i \ \ \text{is a
$G$-torsor on
$\widehat{(X_R)}_{x_i}$, 
}
\\ &\text{$\nu_i$ a trivialization of $\cF_i$ on $\widehat{(X_R)}_{x_i}-x_i$,}
\\
&\text{ and \ $\mu_1$ is a trivialization of
$\cF_1$ on
$\widehat{(X_R)}_{x_2}$}
\endaligned
\right\},
\end{equation}

We again view ${G_{X, \mathcal O}}$  as a group scheme on $X^2$ by pulling it back from the second factor. Then we can define the second action of ${G_{X, \mathcal O}}$  on $\widetilde{\GG_ X\times \GG_ X}$ by letting ${G_{X, \mathcal O}}$ act by altering both of the trivializations $\mu_1$ and $\nu_2$. This action is also free and exhibits $q: \widetilde{\GG_ X\times \GG_ X} \to \GG_ X  \tii \times \GG_ X$ as a ${G_{X, \mathcal O}}$ torsor. Thus, we conclude that the sheaf $\tii \cB$ in \eqref{Global convolution product} exists and is unique.

Let us note that the map
$m$ is a stratified small map -- regardless of the stratification on $X$. To see this,
let us denote by $\Delta \subset X^2$ the diagonal and set $U = X^2-\Delta$.
Then we can take, in  definition \ref{small}, as $W$ the locus of points lying over
$U$. That $m$ is small now follows as
$m$ is an isomorphism over $U$ and over points of $\Delta$ the map $m$
coincides with its analogue in section \ref{Convolution product} which is semi-small
by proposition \ref{stratified semi-small}.

We will now construct the commutativity constraint. For simplicity we specialize
to the case $X=\bA^1$. The advantage is that we can once and for all choose
a global coordinate.  Then the choice of a global coordinate on
$\bA^1$, trivializes $\GG_ X$ over $X$; let us write $\tau :\GG_ X
\to \G$ for the projection. 
By
restricting
$\GG_ {X^{(2)}}$ to the diagonal $\Delta \cong X$ and to $U$, and observing that
these restrictions are isomorphic to $\GG_ X$ and to
$(\GG_ X\times\GG_ X)\rest U$, respectively, we get the following diagram
\begin{equation}
\begin{CD}
\GG_ X   @>i>> \GG_ {X^{2} }@<j<< (\GG_ X\times\GG_ X)\rest U
\\ @VVV @VVV @VVV
\\ X  @>>>  X^2 @<<<\ \ \ U\ \ \ .
\end{CD}
\end{equation}

Let us denote $\tau^o =
\tau^*[1] : \pg \to \operatorname{P}_{G_{X,\mathcal O}}(\GG_ X,\k)$
and
$i^o=i^*[-1]: \operatorname{P}_{G_{X,\mathcal O}}(\GG_ X,\k)\to \pg $.
For $\cA_1,\cA_2\in\pg$ we have:
\begin{equation}\label{convolution products agree}
\aligned &\text{a)}\qquad
\tau^o \cA_1   \ast_X\tau^o \cA_2 \ 
\cong \ 
j_{!*}\left(\phz (
\tau^o\cA_1\Lbten\tau^o\cA_2)
\rest U \right)
\\ &\text{b)}\qquad\tau^o
(\cA_1\ast\cA_2) \ \cong \  i^o 
( \tau^o \cA_1
\ast_X \tau^o\cA_2 ) \,.
\endaligned
\end{equation}

Part a) follows from smallness of $m$ and lemma \ref{small and semi-small maps},
and part b) follows directly from definitions. 

Utilizing the the statements above yields the following sequence of isomorphisms:

\begin{equation}
\gathered
\tau^o(\cA_1\ast\cA_2) \cong i^oj_{!*} \left(
\phz(\tau^o\cA_1\Lbten\tau^o\cA_2)\rest U \right)
\\
\cong  i^oj_{!*}
(\phz(\tau^o\cA_2\Lbten\tau^o\cA_1)\rest U)
\cong\tau^o(\cA_2\ast\cA_1)\,.
\endgathered
\end{equation}

Specializing this isomorphism to (any) point on the diagonal yields a
functorial isomorphism between 
$\cA_1\ast\cA_2$ and $\cA_2\ast\cA_1$. This
gives us a commutativity constraint making $\pg$ into a
tensor category. In the next section we 
modify this commutativity constraint slightly. The modified commutativity constraint will be used in the rest of the
paper.

\begin{rmk}\label{arbitrary curve}
One can avoid having to specialize to the case $X=\bA^1$ here, as well as in the
next section. This can be done, for example, following \cite{BD} and dealing with all choices of a local coordinate at all points of the curve $X$. This gives rise to the
$\operatorname{Aut}(\OO)$-torsor $\hat X \to X$. The functor $\tau^o 
:\pg \to \operatorname{P}_{G_{X,\mathcal O}}(\GG_ X,\k)$ is constructed by
noting that $\GG_X \to X $ is the fibration associated to the 
$\operatorname{Aut}(\OO)$-torsor $\hat X \to X$ and the
$\operatorname{Aut}(\OO)$-action on $\G$. By proposition \ref{another
equivalence}, sheaves in $\pg$ are $\operatorname{Aut}(\OO)$-equivariant and
hence we can transfer them to sheaves on $\GG_X$. 
\end{rmk}

\section{\bf  Tensor functors}
\lab{Tensor functors}

In this section we 
show that our functor
\begin{equation}
{\oh}^* \  \cong \ F \ = \
\bigoplus_{\nu\in X_*(T)} \  \oh^{2\htt(\nu)}_c(S_\nu,-)\ : \  \pg \to \V 
\end{equation}
is a tensor functor. In the case when $\k$ is not a field, the argument is slightly more complicated and we have to make use of some results from section 10. However, the results of this present section are used in section 7 only in the case when $\k$ is a field and not in full generality till chapter 11.

 First, we
alter, following Beilinson and Drinfeld,  the commutativity constraint of the previous section slightly. Let us write $\operatorname{Mod}_\k^\epsilon$ for the tensor
category of finitely generated $\mathbb Z/2\mathbb Z$ -\, graded (super) modules over
$\k$. Let us consider the global cohomology functor as a functor $\mathbb {\oh}^*:
\pg\to \operatorname{Mod}_\k^\epsilon$; here we only keep track of the parity of the
grading on global cohomology.  Then:

\begin{lem}
The functor  ${\oh}^*: \pg\to \operatorname{Mod}_\k^\epsilon$ is a tensor
functor with respect to the commutativity constraint of the previous section. 
\end{lem}

\begin{proof}
We use the interpretation of the convolution product as a fusion product,
explained in the previous section. Let us recall that we write
$\pi:\G_{X^{2}}
\to X^2$ for the projection and again set $X=\bA^1$. The lemma is an immediate
consequence of  the following statements:
\begin{subequations}\label{global tensor functor}
\begin{gather}
R\pi_*(\tau^0(\cA_1)*_X\tau^0(\cA_2))|U  \ \ \text{is the
constant sheaf}\ \
{\oh}^*(\G,\cA_1)\otimes {\oh}^*(\G,\cA_2)\,.
\label{kunneth}
\\ R\pi_*(\tau^0(\cA_1)*_X\tau^0(\cA_2))|\Delta = \tau^0(\mathbb
{\oh}^*(\G,\cA_1*\cA_2))\,.
\label{compatibility}
\\ \text{the sheaves}\ \ \ R^k\pi_*(\tau^0(\cA_1)*_X\tau^0(\cA_2)) \ \ \text{are
constant }\,.\label{constant}
\end{gather}
\end{subequations}
From \eqref{convolution products agree} we immediately conclude \eqref{compatibility} in general and \eqref{kunneth} when $\k$ is field. To prove \eqref{kunneth} in  general, we must show:
\begin{equation}\label{perverse kunneth}
{\oh}^*(\G\times\G,{^p{\oh}^0}(\cA_1\Lbten \cA_2))\ = \ {\oh}^*(\G,\cA_1)\otimes {\oh}^*(\G,\cA_2)\,.
\end{equation}
We will argue this point last and deal with
\eqref{constant} next. Let us write $\tilde \pi : \GG_ X \tii \times \GG_ X \to X^2$ for the natural projection. Then $\tilde \pi = \pi \circ m$. Thus, in order to prove \eqref{constant}  it suffices to show:
\begin{equation}
\text{$R^k\tilde\pi_*\tilde \cB$ is constant;} 
\end{equation}
recall that here $q^*\tilde \cB = p^*(\tau^0(\cA_1)\Lbten\tau^0(\cA_2))$. To do so, we will show that the stratification underlying the sheaf $\tilde \cB$ is smooth over $X^2$. Recall that by a choice of a global coordinate on $X=\bA^1$ we get an isomorphism $\G_X \cong \G \times X$. Thus, the sheaves $\tau^0(\cA_1)$ and $\tau^0(\cA_2)$ are constructible with respect to the stratification $\G^\la_X$ which correspond to $\G^\la \times X$ under the above isomorphism; here, as usual, $\lambda\in X_*(T)$. These strata are smooth over the base $X$ by construction. Thus, we conclude that the sheaf $\tilde\cB$ is constructible with respect to the strata $\G^\la_X \tii\times \G^\mu_X$, for $\lambda,\mu\in X_*(T)$, which are uniquely described by the following property:
\begin{equation}
q^{-1}(\G^\la_X \tii\times \G^\mu_X) \ = \ p^{-1}(\G^\la_X\times \G^\mu_X)\,.
\end{equation}
In other words, the strata  $\G^\la_X \tii\times \G^\mu_X$ are quotients of $p^{-1}(\G^\la_X\times \G^\mu_X)$ by the second ${G_{X, \mathcal O}}$ action on $\widetilde{\GG_ X\times \GG_ X}$ defined in section \ref{The commutativity constraint} which makes  $q: \widetilde{\GG_ X\times \GG_ X} \to \GG_ X  \tii \times \GG_ X$  a ${G_{X, \mathcal O}}$ torsor. As such, the $\G^\la_X \tii\times \G^\mu_X$ are smooth. Furthermore, the projection morphism $\tilde \pi_{\la,\mu} : \G^\la_X \tii\times \G^\mu_X \to X^2$ is smooth. This can be verified either by a direct inspection or concluded by general principles from the fact that all the fibers of $\tilde \pi_{\la,\mu}$ are smooth and equidimensional. This, then, lets us conclude \eqref{constant}. 

It remains to argue \eqref{perverse kunneth}. Let us first assume that one of the factors ${\oh}^*(\G,\cA_i)$ is flat over $\k$. Then, by Lemma  \eqref{outer tensor product is perverse}, the sheaf $\cA_1\Lbten \cA_2$ is perverse. Then, again using the flatness of ${\oh}^*(\G,\cA_i)$, we get
\begin{multline}
{\oh}^*(\G\times\G,{^p{\oh}^0}(\cA_1\Lbten \cA_2))\ = \ {\oh}^*(\G\times\G,\cA_1\Lbten \cA_2)\ = 
\\
{\oh}^*(\G,\cA_1)\overset{L}\otimes {\oh}^*(\G,\cA_2)\ = \ {\oh}^*(\G,\cA_1)\otimes {\oh}^*(\G,\cA_2) \,.
\end{multline}

To argue the general case we make use of Corollary  \ref{enough projectives} and Proposition \ref{structureProjectives}.   Corollary  \ref{enough projectives} allows us to write any $\cA\in\pg$ as a quotient of a projective $\cP\in\pp Z$ and Proposition \ref{structureProjectives} tells us that ${\oh}^*(\G,\cP)$ is free over $\k$; here $Z$ is any $\GO$-invariant finite dimensional subvariety of $\G$ which contains the support of $\cA$. Let us consider a resolution of $\cA_1$ by such projectives:
\begin{equation}
\cQ \to \cP \to \cA_1 \to 0 \,.
\end{equation}
As the functor $\cA \mapsto {^p{\oh}^0}(\cA\Lbten \cA_2)$ is right exact, we get an exact sequence
\begin{equation}
{^p{\oh}^0}(\cQ\Lbten \cA_2)\to {^p{\oh}^0}(\cP\Lbten \cA_2)\to {^p{\oh}^0}(\cA_1\Lbten \cA_2) \to 0\,.
\end{equation}
Because cohomology is a an exact functor and making use of the fact that we have already proved \eqref{perverse kunneth} for the first two terms, we get an exact sequence
\begin{multline}
\ \ \ \ \ \ \ \ \ {\oh}^*(\G,\cQ)\ten{\oh}^*(\G,\cA_2)\to {\oh}^*(\G,\cP)\ten{\oh}^*(\G,\cA_2)\to
\\
\to{\oh}^*(\G\times\G, {^p{\oh}^0}(\cA_1\Lbten \cA_2))\to 0\,.
\end{multline}
Comparing this exact sequence to the one we get by tensoring the exact sequence
\begin{equation}
 {\oh}^*(\G,\cQ) \to  {\oh}^*(\G,\cP) \to  {\oh}^*(\G,\cA_1) \to 0
\end{equation}
with ${\oh}^*(\G,\cA_2)$ concludes the proof.
\end{proof}

\begin{rmk}
The statements in \eqref{global tensor functor} hold for an arbitrary
curve $X$. This can be seen by utilizing the $\operatorname{Aut}(\OO)$-torsor
$\hat X \to X$ of remark \ref{arbitrary curve} and proposition \ref{another
equivalence}; for details see \cite{Na}.
\end{rmk}

Let $\V$ denote the category of finite dimensional  vector spaces over
$\k$. To make ${\oh}^*: \pg \to \V$ into a tensor functor we
alter, following Beilinson and Drinfeld,  the commutativity constraint of the
previous section slightly.
We consider the constraint from section \ref{The commutativity constraint}  on the category
$\pg\otimes\operatorname{Mod}_\k^\epsilon$ and restrict it
to a  subcategory that we
identify with
$\pg$.
Divide $\GG$ into unions of connected
components $\GG=\ \GG_+\cup\GG_-$ so that the dimension of  $\GO$-orbits
is even in $\GG_+$ and odd in $\GG_-$. This gives a $\Z_2$-grading
on the category   $\pg$ hence a new $\Z_2$-grading
on $\pg\otimes\operatorname{Mod}_\k^\epsilon$. The subcategory of even
objects is identified with $\pg$ by forgetting the grading. Hence, we conclude from
the previous lemma:

\begin{prop} The functor ${\oh}^*: \pg \to
\V$ is a tensor functor with respect to the above
commutativity constraint.
\end{prop}

Let us write $\V(X_*(T)))$ for the (tensor) category of finitely generated
$\k$-modules with a $X_*(T)$-grading. We can view $F = \oplus_{\nu\in
X_*(T)} \  F_\nu$ as a functor from $ \pg$ to $\V(X_*(T)))$. Then we have the
following generalization of the previous proposition:

\begin{prop}  \label{weighted tensor functor}
The functor $F: \pg \to
\V(X_*(T)))$ is a tensor functor.
\end{prop}

\begin{proof}
The notion of the subspaces $S_\nu$ and $T_\nu$ can be extended to the situation of
families, i.e., to the global Grassmannians $\GG_ {X^{n}}$. Recall that the fiber of the
projection $r_n:\GG_ {X^{n}} \to X^n$ over the point
$(x_1,\dots,x_n)$ is simply $\prod_{i=1}^k \GG_ {y_i}$\,, where
$\{y_1,\dots,y_k\}=\{x_1,\dots,x_n\}$, with all the $y_i$ distinct. Attached to the
coweight $\nu\in X_*(T)$ we associate the ind-subscheme
\begin{equation}\label{products of Snu}
 \prod_{\nu_1+\dots+\nu_k=\nu} S_{\nu_i} \ \subset \ \prod_{i=1}^k \GG_ {y_i}=
r_n^{-1}(x_1,\dots,x_n)
\end{equation}
These ind-schemes altogether form an ind-subscheme $S_\nu(X^n)$ of $\GG_
{X^{n}}$. This is easy to see for $n=1$ by choosing a global parameter, for example.
By the same argument we see that outside of the diagonals $S_\nu(X^n)$ form a
subscheme. It is now not difficult to check that the closure of this locus lies inside
$S_\nu(X^n)$. Similarly, we define the ind-subschemes $T_\nu(X^n)$. Let us write
$s_\nu$ and $t_\nu$ for the inclusion maps of $S_\nu(X^n)$ and $T_\nu(X^n)$ to
$\GG_ {X^{n}}$, respectively. We have the action of $\bG_m$ on
$\GG_ {X^{n}}$ via the cocharacter $2\check\rho$. The fixed point set of this action
consists of the locus of products of the fixed points in the individual affine
Grassmannians, i.e., above the point $(x_1,\dots,x_n)$ where
$\{x_1,\dots,x_n\}=\{y_1,\dots,y_k\}$, with all the $y_i$ distinct, the fixed points are of
of the form
\begin{equation}
(L_{\nu_1}, \dots, L_{\nu_k}) \in \prod_{i=1}^k \GG_ {y_i};
\end{equation}
recall that we write $L_\nu$ for the point in $\G$ corresponding to the cocharacter
$\nu\in X_*(T)$. We write $C_\nu$ for the  subset of the fixed point locus lying inside
$S_\nu(X^n)$, i.e., 
\begin{equation}
C_\nu\cap  r_n^{-1}(x_1,\dots,x_n) \ = \ \bigcup_{\nu_1+\dots+\nu+k=\nu}
\{(L_{\nu_1}, \dots, L_{\nu_k})\}\,.
\end{equation}
Let us write $i_\nu: S_\nu(X^n) \rightarrow \GG_
{X^{n}}$ and $k_\nu: T_\nu(X^n) \rightarrow \GG_
{X^{n}}$ for the inclusions. By the same argument as in the proof of theorem
\ref{dimension estimates} we see  that 
\begin{align}
S_\nu(X^n) \ &= \ \{z\in\GG_{X^{n}}\mid \lim_{s\to 0} 2\check\rho(s)z \in C_\nu \}
\\
\intertext{and}
T_\nu(X^n) \ &= \ \{z\in\GG_{X^{n}}\mid \lim_{s\to \infty} 2\check\rho(s)z \in C_\nu
\}\,.
\end{align}
Let us write $p_\nu: S_\nu(X^n) \to C_\nu$ and $q_\nu: T_\nu(X^n) \to C_\nu$ for the
retractions:
\begin{align}
p_\nu(z) \ &= \ \lim_{s\to 0} 2\check\rho(s)z \qquad \text{for} \ \ z \in S_\nu(X^n)
\\
q_\nu(z) \ &= \ \lim_{s\to \infty} 2\check\rho(s)z \qquad \text{for} \ \ z \in T_\nu(X^n)\,.
\end{align}
 Furthermore,
\begin{equation}
C_\nu \ = \ S_\nu(X^n) \cap T_\nu(X^n)\,.
\end{equation}
By  Theorem 1 of \cite{Br} we conclude that 
\begin{equation}
i_\nu^!s_\nu^*\cB \ = \ k_\nu^* t_\nu^!\cB \qquad \text{for} \ \
\cB\in\operatorname{P}_{G_{X^n, \mathcal O}}(\GG_ {X^n},\k)\,.
\end{equation}
Rephrasing this result in terms of the contractions $p_\nu$ and $q_\nu$ gives:
\begin{equation}\label{global updown}
R(p_\nu)_!s_\nu^*\cB \ = \ R(q_\nu)_* t_\nu^!\cB \qquad \text{for} \ \
\cB\in\operatorname{P}_{G_{X^n, \mathcal O}}(\GG_ {X^n},\k)\,.
\end{equation}

Let us now, for simplicity, choose $X=\bA^1$. Let $\cA_1,\cA_2\in\pg$. We write
$\cB_1= \tau^o\cA_1$ and $\cB_2= \tau^o\cA_2$ and form the convolution product
$\cB_1\ast\cB_2 = Rm_*\tii\cB$.  By statement
\eqref{constant} we see:
\begin{equation}\label{constant2}
\begin{split}
\text{The sheaf} \ \  R^k\pi_*Rm_*\tii\cB \ \ \ \text{on $X^2$ is constant with fiber} 
\\
 {\oh}^k(\G,\cA_1\ast\cA_2) \ = \ \bigoplus_{k_1+k_2=k}{\oh}^{k_1}(\G,\cA_1)
\otimes {\oh}^{k_2}(\G,\cA_2) \,.
\end{split}
\end{equation}

Let us now consider the sheaves 
\begin{equation}
\cL^k_\nu(\cA_1,\cA_2)=R^k\pi_*R(p_\nu)_!s_\nu^*Rm_*\tii\cB = 
R^k\pi_*R(q_\nu)_*t_\nu^!Rm_*\tii\cB\,;
\end{equation}
here we have used \eqref{global updown} to identify the last two sheaves. Let us calculate the stalks of this sheaf. First, by definition, 
\begin{equation}
\cL_\nu^{k}(\cA_1,\cA_2)_{(x_1,x_2)} = 
\begin{cases} 
{\oh}^{k}_c(S_\nu, \cA_1\ast\cA_2) &\text{if $x_1=x_2$}
\\
\oplus_{\nu_1+\nu_2 = \nu}{\oh}^{k}_c(S_{\nu_1}\times S_{\nu_2},{^p{\oh}^0}(\cA_1\Lbten \cA_2) )
&\text{if $x_1\neq x_2$}\,.
\end{cases}
\end{equation}
Arguing in the same way as in the proof of \eqref{perverse kunneth} we see that 
\begin{equation}
{\oh}^*_c(S_{\nu_1}\times S_{\nu_2},{^p{\oh}^0}(\cA_1\Lbten \cA_2))\ = \ {\oh}^*_c(S_{\nu_1},\cA_1)\otimes {\oh}^*_c(S_{\nu_2},\cA_2)\,.
\end{equation}

We conclude that
\begin{subequations}\label{stalks}
\begin{equation}
\cL^k_\nu(\cA_1,\cA_2)_{(x_1,x_2)} \ = \ 
0 \qquad \text{if $k\neq 2\rho(\nu)$}\,,
\end{equation}
and
\begin{multline}
\cL^{2\rho(\nu)}_\nu(\cA_1,\cA_2)_{(x_1,x_2)} = 
\\
\begin{cases} 
{\oh}^{2\rho(\nu)}_c(S_\nu, \cA_1\ast\cA_2) &\text{if $x_1=x_2$}
\\
\oplus_{\nu_1+\nu_2 = \nu}{\oh}^{2\rho(\nu_1)}_c(S_{\nu_1},
\cA_1)\otimes{\oh}^{2\rho(\nu_2)}_c(S_{\nu_2}, \cA_2)&\text{if $x_1\neq x_2$}\,.
\end{cases}
\end{multline}
\end{subequations}

We now proceed as in the proof of theorem \ref{exactness}. Let us consider the
closures $\barr{ S_\nu}(X^n)$ and $\barr{ T_\nu}(X^n)$ 
of the ind-subschemes $S_\nu(X^n)$
and $T_\nu(X^n)$ and let us write $\barr{ i_\nu}: 
\barr{ S_\nu}(X^n) \rightarrow \GG_
{X^{n}}$ and $\barr {k_\nu}: \barr{ T_\nu}(X^n) \rightarrow \GG_
{X^{n}}$ for the inclusions. 
Let us write $\cB=Rm_*\tii\cB$. Then we have the
following canonical morphisms
\begin{subequations}
\begin{align}
R\pi_*\cB=R\pi_!\cB &\rightarrow R\pi_!\barr{ s_\nu}^*\cB \label{s map}
\\
R\pi_*\cB &\leftarrow R\pi_*t_\nu^!\cB \label{t map}\,.
\end{align}
\end{subequations}
These morphisms give us two filtrations of
$R^k\pi_*Rm_*\tii\cB$, one by kernels of the the morphisms
\begin{equation}
R^k\pi_*\cB \rightarrow R^k\pi_!\barr{s_\nu}^*\cB
\end{equation}
 and the other by images of the morphisms
\begin{equation}
R^k\pi_*\barr{ t_\nu}^!\cB \rightarrow
R^k\pi_*\cB\,.
\end{equation}
By the discussion above, these filtrations are complementary and hence yield the
following canonical isomorphism
\begin{equation}
R^k\pi_*Rm_*\tii\cB \ = \ \bigoplus_{2\rho(\nu)=k} \cL^k_\nu(\cA_1,\cA_2)\,.
\end{equation}

By \eqref{constant2} the sheaf on the left hand side is constant. Therefore the sheaves $
\cL^k_\nu(\cA_1,\cA_2)$ must be constant. Appealing to  \eqref{stalks} completes the
proof. 

\end{proof}

\section{\bf The case of a field of characteristic zero}
\label{The case of a field of characteristic zero}

In this section we treat the case when the base ring $\k$ is a field of characteristic
zero. This case was treated already in \cite{Gi} when $\k=\bC$. Here we make use of Tannakian
formalism, using \cite{DM} as a general reference. In section \ref{Construction
of the group scheme}, where we work over an arbitrary base ring $\k$, we carry
out the constructions explicitly without referring to the general Tannakian formalism.

\begin{lem}\label{semisimplicity}
If  \,$\k$ is a field of characteristic zero then the category $\pg$ is semisimple. In
particular, the sheaves $\II_!(\la,\k)$, $\II_{!*}(\la,\k)$,
and $ \II_*(\la,\k)$ are isomorphic. 
\end{lem}

\begin{proof}
The parity vanishing of the stalks of  $\II_{!*}(\la,\k)$, proved in 
\cite{Lu}, section
11, and the fact that the orbits $\GG^\la$ are simply connected implies immediately
that there are no extensions between the simple objects in
$\ps (\G,\k)$. Thus, there are no extensions in the full subcategory $\pg$ (which then, obviously, coincides with $\ps (\G,\k)$). 
\end{proof}

\begin{rmk}
The use of the above lemma can be avoided. One must then ignore this section and first go through the rest
of the paper in the case when $\k$ is a field of characteristic zero. The arguments of section
\ref{identification}, in a greatly simplified form, then give theorem \ref{characteristic zero}. 
\end{rmk}

The constructions above  and the properties we established
suffice for
verifying the conditions of
the proposition 1.20 in
\cite{DM} and then also the conditions of the  theorem 2.11 in \cite{DM},
which are summarized by the phrase\ 
``$(\pg,\ast,F)$ is a  {\em neutralized Tannakian category}''.
Hence, by theorem 2.11 of \cite{DM}, we conclude:
\begin{equation}
\begin{gathered}
\text{there is a group scheme $\tii G$ over $\k$ such that}
\\
\text{the
category of finite dimensional $\k$-representations of $\tii G$}
\\
\text{is equivalent to
$\pg$, as tensor categories}\,.
\end{gathered}
\end{equation}

We will now identify the group $\tii G$. Let us write $\check G$ for the dual
group of $G$, i.e., $\check G$ is the split reductive group over $\k$ whose root
datum is dual to that of $G$.

\begin{thm} \label{characteristic zero}
The category of finite dimensional $\k$-representations of $\check G$ is
equivalent  to $\pg$, as tensor categories.
\end{thm}

Before giving a proof of this theorem we discuss it briefly from the point of view of representation theory.
We can view the theorem as giving us a geometric interpretation of representation theory of $\check G$.
First of all, as we use global cohomology as fiber functor, it follows that the representation space for the
representation $V_\cF$, associated to $\cF\in\pg$,  is the global cohomology $\oh^*(\G,\cF)$.
As the proof of the theorem will show, $\check T$, the dual torus of $T$, is a maximal torus in $\check G$.
The decomposition of  $V_\cF$ into its $\check T$-weight spaces is then given by theorem \ref{exactness}:
\begin{equation}
\mathbb {\oh}^*(\G,\cF) \ 
\cong  \ \bigoplus_{\nu\in X^*(\check T)} \ 
\oh^{2\htt(\nu)}_c(S_\nu,\cF)\,.
\end{equation}
Given a dominant
$\lambda\in X_*(T) = X^*(\check T)$ we can associate to it both the highest weight representation $L(\la)$
of $\check G$ and the sheaf $\II_{!*}(\la,\k)\in\pg$.  Obviously,
$V_{\II_{!*}(\la,\k)}$ is irreducible and by the formula above we see that it is of highest weight $\la$.
Hence, $V_{\II_{!*}(\la,\k)} = L(\la)$. Combining this discussion with lemma \ref{semisimplicity} and
proposition \ref{CanonicalBasis} gives:

\begin{cor}
The $\nu$-weight space $L(\la)_\nu$ of $L(\la)$ can be canonically identified 
with the $\k$-vector space
spanned by the irreducible components of ${\barr{\GG_\la}}\cap S_\nu$.
In particular, the dimension of
$L(\la)_\nu$ is given by the 
number of irreducible components of ${\barr{\GG_\la}}\cap S_\nu$\,.
\end{cor}

The rest of this section is devoted to the proof of theorem \ref{characteristic zero}. We begin with an
observation:
\begin{equation}
\text{the group scheme $\tii G$ is a split connected reductive algebraic group}
\end{equation}

To see that  $\tii G$ is algebraic, we observe that it has a tensor generator. Let
$\la_1, \dots, \la _r$ be a set of generators for the dominant weights in $X_*(T)$. As a
generator we can then take $\oplus \II_{!*}(\la_i,\k)$. It is tensor generator because for
any dominant $\la$ the sheaf $\II_{!*}(\la,\k)$ appears as a direct summand in the
product
\begin{equation}
\II_{!*}(\la_1,\k)^{\ast
k_1} \ast \dots \ast\II_{!*}(\la_r,\k)^{\ast k_r}\,;
\end{equation}
here $\la = \sum k_1\la_1 + \dots + k_r\la_r$. Thus, by \cite{DM}, proposition 2.20,
$\tii G$ is an algebraic group. As there is no  tensor
subcategory of $\pg$ whose objects are direct sums of finitely many fixed
irreducible objects the group $\tii G$ is connected by \cite{DM}, corollary 2.22.
Finally, as $\pg$ is semisimple, $\tii G$ is reductive, by \cite{DM}, proposition
2.23. To see that $\tii G$ is split, we exhibit a split maximal torus in $\tii G$.
By proposition \ref{weighted tensor functor} the fiber functor $F ={\oh}^*$
factors as follows:
\begin{equation}
F ={\oh}^*: \pg \to \V(X_*(T))) \to \V\,.
\end{equation}
This gives us a homomorphism $\check T \to \tii G$; here $\check T$ is the
torus dual to $T$. As any character $\la\in X^*(\check T) = X_*(T)$ appears as the 
direct summand $F_\lambda(\II_{!*}(\la,\k))$ in $F(\II_{!*}(\la,\k))$ we conclude
that $\check T$ is a split torus in $\tii G$. It is clearly maximal as the
representation ring of $\tii G$ is of the same rank as $\check T$. 

It now remains to identify the root datum of $\tii G$ with the dual of the root
datum of $G$. Recall that we have also fixed a choice of positive roots,
i.e., a Borel $B$ such that $T\subset B\subset G$. The root datum of $G$ is then
given as $(X^*(T), X_*(T), \De(G,T), \check \De(G,T))$, 
where $\De(G,T)
\subset X^*(T)$ are the roots and 
$\check\De(G,T) \subset X_*(T)$ are the
coroots of $G$ with respect to $T$. Because
$X^*(\check T) = X_*(T)$ and $X_*(\check T) = X^*(T)$, it suffices to show that
\begin{equation}\label{equality of root data}
\text{
$\De(\tii G, \check T) = \check \De(G,T)$ 
\ \ \ and \ \ \ $\check \De(\tii G,
\check T) =
\De(G,T)$}\,. 
\end{equation}
To this end we note that theorem \ref{dimension estimates},
corollary \ref{weight picture}, and proposition \ref{CanonicalBasis} imply that:

\begin{equation}\label{a1}
\begin{gathered}
\text{The irreducible representations of $\tii G$ are}
\\
\text{parameterized by dominant
coweights $\la\in X_*(T)$.}
\end{gathered}
\end{equation}
and
\begin{equation}\label{a2}
\begin{gathered}
\text{The $\check T$-weights of the irreducible
representation $L(\lambda)$}
\\
\text{associated to $\la$ are the same as the $\check T$-weights}
\\
\text{of the irreducible
representation of $\check G$ associated to $\la$}\,.
\end{gathered}
\end{equation}

We now argue using the pattern of the weights. For clarity we spell out this familiar structure. From \eqref{a2} we conclude:
 
\begin{subequations}\label{structure of weights in characteristic zero}
\begin{equation}
\text{The weights of $L(\lambda)$ are symmetric under the Weyl group $W$.}
\end{equation}
\begin{equation}
\begin{gathered}
\text{For $\alpha$ a simple positive root of $G$, with the corresponding coroot $\check\alpha$, the}
\\
\text{weights of $L(\lambda)$ on the line segment between $s_{\check \alpha}(\lambda)$ and $\lambda$ are the weights}
\\
\text{on that segment which are of the form $\lambda-k\check\alpha$, with $k$ an integer.}
\end{gathered}
\end{equation}
\end{subequations}

Note that the choice of a Borel subgroup of $\tii G$ is equivalent to a consistent
choice of a line, the highest weight line,  in each irreducible representation of $\tii
G$. The choice  $F_\la(\II_{!*}(\la,\k))$ in $F(\II_{!*}(\la,\k))$ for all
dominant $\la\in X_*(T)$ yields a Borel subgroup $\tii B$ of $\tii G$ such that
the dominant weights of $\tii G$ in $X^*(\check T)$ coincide with the dominant
coweights of $G$ in $X_*(T)$. This implies that the simple coroot directions
of the triple $(\check T, \tii B,\tii G)$ coincide with the simple root
directions of $(T, B, G)$. The statements \eqref{structure of weights in characteristic zero} above now imply
that the simple roots of the triple $(\check T, \tii B, \tii G)$ coincide with the
simple coroots of
$(T, B, G)$. This, finally gives \eqref{equality of root data}.

\section{\bf Standard sheaves}
\lab{Canonical basis}

In this section we prove some basic results about standard sheaves which will be
crucial for us later. Let us write $\D$ for the Verdier duality functor. 

\begin{prop}
\label{tenk}
We have
\begin{align}
\tag{a}\II_!(\la,\k)\ &\cong\ \II_!(\la,\Z) \bb{\Z}{\aa{L}\ten}\k
\\
\tag{b}\II_*(\la,\k)\ &\cong\ \II_*(\la,\Z) \bb{\Z} {\aa{L} \ten} \k
\\
\tag{c}\D\, \II_!(\la,\k)\ &\cong \II_*(\la,\k)\,.
\end{align}
\end{prop}

\begin{proof}
The proofs of (a) and (b) are analogous and hence we will only prove (a). Because
\begin{equation}
{\oh}^*_c(S_\nu,\II_!(\la,\Z)\bb{\Z}{\aa{L}\ten}\k)
\ =\
{\oh}^*_c(S_\nu,\II_!(\la,\Z)) \bb{\Z} {\aa{L}\ten}\k
\end{equation}
and,  by proposition \ref{CanonicalBasis}, ${\oh}^*_c(S_\nu,\II_!(\la,\Z))$  is a free
abelian group in degree
$2\htt(\nu)$  we conclude that ${\oh}^*_c(S_\nu,\II_!(\la,\Z)\bb{\Z}{\aa{L}\ten}\k)$ is
nonzero only in degree $2\htt(\nu)$. Hence, by lemma \ref{perversity criterion}, we see 
that
$
\II_!(\la,\Z)
\bb{\Z}{\aa{L}\ten}\k
$ is perverse. There is a canonical map
\begin{equation}
\II_!(\la,\k) \xrightarrow{\ \io\ }
\II_!(\la,\Z)\bb{\Z}{\aa{L}\ten}\k\,,
\end{equation}
which is an isomorphism when restricted to $\G^\la$. Therefore, 
applying the  functor $F_\nu$ to the morphism $\io$ and using the 
proposition
\ref{CanonicalBasis}
yields an isomorphism
\begin{equation}
F_\nu(\II_!(\la,\k)) 
=\ 
\k[Irr(\barr {\GG_\la}\cap S_\nu)] 
\xrightarrow{\ F_\nu(\io)\ }
\Z[Irr( {\barr{\GG_\la}} \cap S_\nu)]
\bb{\Z}\ten \k
\ =\
F_\nu(\II_!(\la,\Z)\bb{\Z}{\aa{L}\ten}\k)\,.
\end{equation}
By corollary \ref{faithful and exact} the functor $F= \oplus F_\nu$ is faithful and thus 
we conclude that $\io$ is an isomorphism.

The proof of part (c) proceeds in a similar fashion. First we observe that 
\begin{equation}
{\oh}^*_{T_\nu}(\G,\D\ \II_!(\la,\k))\cong \
\D({\oh}^*_c(T_\nu,\II_!(\la,\k)))
=\
\D({\oh}^*_{c}(S_{w_0\cd\nu},\II_!(\la,\k)))\,.
\end{equation}
Because ${\oh}^*_{c}(S_{w_0\cd\nu},\II_!(\la,\k))$ is a free $\k$-module
concentrated in degree $2\htt(w_0\cd\nu)$, we conclude that
$\D({\oh}^*_{T_{\nu}}(\GG,\II_!(\la,\k)))$ 
is a concentrated in degree
$-2\htt(w_0\cd\nu)= 2\rho(\nu)$. Thus, we conclude that $\D\,\II_!(\la,\k)$ is perverse.
Furthermore, we note that
\begin{equation}\label{intermediate}
\begin{CD}
{\oh}^*_{T_\nu}(\G,\D\ \II_!(\la,\k)) @>\cong>>
\D({\oh}^*_c(T_\nu,\II_!(\la,\k)))
\\
@VVV @VV{\cong}V
\\
{\oh}^*_{T_\nu}(\G^\la,\D\ \II_!(\la,\k)) @>\cong>>
\D({\oh}^*_c(T_\nu\cap\G^\la,\II_!(\la,\k)))\,,
\end{CD}
\end{equation}
which implies that the left hand arrow is also an isomorphism.
We have a canonical map
\begin{equation}
\D\ \II_!(\la,\k)\
\xrightarrow{\ \io\ }
\II_*(\la,\k)
\end{equation}
To show that this map is an isomorphism it suffice to show that  the maps $F_\nu(\io)$
are isomorphisms. Restricting to $\G^\la$ gives us the following commutative diagram:
\begin{equation}
\begin{CD}
{\oh}^*_{T_\nu}(\G,\D\ \II_!(\la,\k)) @>{F_\nu(\io)}>>
{\oh}^*_{T_\nu}(\G,\II_*(\la,\k))
\\
@V{\cong}VV @VV{\cong}V
\\
{\oh}^*_{T_\nu}(\G^\la,\D\ \II_!(\la,\k)) @>\cong>>
{\oh}^*_{T_\nu}(\G^\la,\II_*(\la,\k))\,.
\end{CD}
\end{equation}
In this diagram the bottom arrow is an isomorphism because $\io$ restricted to $\G^\la$ is
an isomorphism, the left vertical arrow is an isomorphism by \eqref{intermediate}, and
finally, the right vertical arrow is an isomorphism by proposition \ref{CanonicalBasis}
(or, rather, by the proof thereof). This shows that $F_\nu(\io)$ is an isomorphism.

\end{proof}

\begin{prop}
\label{!=!*}
The canonical map 
$
\II_!(\la,\Z)
\ra\
\II_{!*}(\la,\Z)
$\
is an isomorphism.
\end{prop}

\begin{proof}
Let us consider the following commutative diagram:
\begin{equation}
\begin{CD}
\II_!(\la,\Z)
@>{\al}>>
\II_{*}(\la,\Z)
\\
@VVV
@VVV
\\
\II_!(\la,\Q)
@>{}>>
\II_{*}(\la,\Q)\,.
\end{CD}
\end{equation}
The bottom map is an isomorphism by lemma \ref{semisimplicity}.
Let us apply the functor $F_\nu$ to this diagram. 
The columns become inclusions by
proposition
\ref{tenk} and the bottom arrow is an isomorphism as we just observed.
Therefore, $F_\nu(\al)$ is an inclusion and 
then so is $\al$.
This implies that  the canonical surjection
$
\II_!(\la,\Z)
\sur\
\II_{!*}(\la,\Z)
$ is an isomorphism.
\end{proof}

\section{\bf Representability of the weight functors}
\lab{Representability of the weight functors}
In section \S\ref{Semi-infiniteOrbits}
 we showed that the functors $F_\nu \df\
\oh^{2\htt(\nu)}_c(S_\nu,\ \,)
\cong\
\oh_{T_\nu}^*(\G,-)
[2\htt(\nu)]
$,
from $ \pg$ to $\V$, for $\nu\in X_*(T)$,  are exact. 
Hence, one would expect them to be (pro)
representable. Here we prove that this is indeed the case:

\begin{prop} 
\label{representable}
Let $Z\subset\G$ be a closed subset which is a finite union of $\GO$-orbits.
The functor $F_\nu$ restricted to $\pp Z$ is represented by a projective
object $\PPP_Z(\nu,\k)$ of $\pp Z$.
\end{prop}

\begin{proof}
We make use of the induction functors. Let us recall their construction. For more details see,
for example, \cite{MiV1}. Let $A$ be an algebraic group acting on a variety $Y$ and let $B$ be
a subgroup of $A$. The forgetful functor $\FF^A_B:\ D_A(Y,\k)\ra\ D_B(Y,\k)$ has a left
adjoint
$\ga^A_B:\ D_B(Y,\k)\ra\ D_A(Y,\k)$ which can be constructed as follows. Consider the
diagram
\begin{equation}
Z \xleftarrow{ \ \ p \ \ } A\times Z \xrightarrow{\ \ q\ \ } A\times_B Z \xrightarrow{\ \ a\ \ }
Z\,. 
\end{equation}
The maps $p$ and $q$ are projections, and $a$ is the action map. The group $A\times B$ acts
on the leftmost copy of $Z$ via the factor $B$, on $A\times Z$ and $A\times_B Z$ by the
formula $(a,b)\cdd(a',z)=(a\cdd a' \cdd b\inv,b\cdd z)$, and on the leftmost copy of $Z$ via
the factor $A$. The left adjoint $\ga^A_B$ is now given by 
\begin{equation}
\ga^A_B(\AA) \ = \ \ a_!\ti \AA \ \ \text{where} \ \ \ \ti \AA \ \ \text{is defined via}\ \ q^!\ti
\AA\ \cong\ p^!\AA\,; \text{for} \ \ \AA\in D_B(Y,\k)
\end{equation}

Let $\cO_n=\cO/z^{n+1}$ and let us write $G_{\cO_n}$ for the 
algebraic group whose
$\bC$-points are $G(\cO_n)$. We use analogous notation for other 
groups. Now choose
$n>>0$ so that the $\GO$-action on $Z$ factors through the action of
$G_{\cO_n}$. We write $\PPP_Z(\nu,\k) =
{}^p\oh^0(\gamma_{\{e\}}^{G_{\cO_n}}(\k_{T_\nu\cap Z}[-2\htt(\nu)] ))$. We claim: 

\begin{equation}
\label{representability of Fsubnu}
\text{the functor $F_\nu: \pp Z\to \V$ is represented by $\PPP_Z(\nu,\k)$}
\end{equation}

To see this, let $\cA\in \pp Z$, and then:

\begin{multline}
\label{representability}
F_\nu(\mathcal A) \ =\ {\oh}^{2\htt(\nu)}_{T_\nu}(Z,\cA)  \ =\ 
\\
\Ext_{\oD(Z,\k)}^0(\k_{T_\nu\cap
Z}[-2\htt(\nu)],
\cF^{G_{\cO_n}}_{\{e\}}\ \cA)\  \cong
\\
\Ext_{\oD_{G_{\cO_n}}(Z,\k)}^0(\gamma^{G_{\cO_n}}_{\{e\}}\k_{T_\nu\cap
Z}[-2\htt(\nu)],
\ \cA)\,.
\end{multline}

Let us write
$\cF=\gamma^{G_{\cO_n}}_{\{e\}}\k_{T_\nu\cap Z}[-2\htt(\nu)]$.
Then
\begin{equation}
\label{leftexact}
\text{the sheaf}  \ \ \cF \ \ \text{lies in} \ \ ^p\oD^{\le 
0}_{G_{\cO_n}}(Z,\k)\,.
\end{equation}
To see this, let us write $d$ for the largest integer such that $\ 
^p\oH^d(\cF)\ne 0$.
Then, we see, as in \eqref{representability}, that  
\
\begin{multline}
0\ne\Hom_{\oD_{G_{\cO_n}}(Z,\k)}(\cF,{}^p\oh^d(\cF)[-d])
=
\\
\Ext^{-d}_{\oD_{G_{\cO_n}}(Z,\k)}(\cF,
{}^p\oh^d(\cF))\cong{\oh}^{2\htt(\nu)-d}_{T_\nu}(Z,{}^p\oh^d(\cF))\,.
\end{multline}
According to theorem
\ref{weightfunctors}
this forces $d=0$ and we have proved \eqref{leftexact}. This immediately implies
\eqref{representability of Fsubnu}.

\end{proof}

\begin{cor} \label{enough projectives}
The category
$\oP_\GO(Z,\k)$ has enough projectives.
\end{cor}

\begin{proof}
Let $\cA\in
\oP_\GO(Z,\k)$.
Choose finitely generated $\Bbbk$-projective covers
$f_\nu:P_\nu\ra F_\nu(\cA)$.
Then
\begin{equation}
\Hom(
P_\nu\otimes_\Bbbk
\PPP_Z(\nu,\k),
\cA)
\
\cong
\
\Hom_\k[
P_\nu,
\Hom(\PPP_Z(\nu,\k),
\cA)]
\
\cong
\
\Hom_\k[
P_\nu,
F_\nu(
\cA)]\,.
\end{equation}
By construction, the map $p_\nu\in\Hom(P_\nu\otimes_\Bbbk\PPP_Z(\nu,\k),\cA)$
corresponding to the  $\Bbbk$-projective cover $f_\nu:P_\nu\ra F_\nu(\cA)$ satisfies
$F_\nu(p_\nu)=f_\nu$. 
Since $\pl_\nu F_\nu=F$ is exact and faithful,
$\oplus_\nu\ P_\nu\otimes_\Bbbk \PPP_Z(\nu,\k) $ is a projective cover of $\cA$.

\end{proof}

We can describe the sheaf
$\PPP_Z(\nu,\k)={}^p\oh^0(\gamma_{\{e\}}^{G_{\cO_n}}(\k_{T_\nu\cap Z}[-2\htt(\nu)] ))$
rather explicitly as follows. Let us consider the following diagram:
\begin{equation}
\begin{CD}
Z @<p<< G_{\cO_n} \times Z @>a>> Z
\\
@AiAA @AAA @|
\\
T_\nu\cap Z @<r<<   G_{\cO_n} \times (T_\nu\cap Z) @>{\tii{a}}>> Z
\end{CD}
\end{equation}
and use it to calculate $\gamma_{\{e\}}^{G_{\cO_n}}(\k_{T_\nu\cap Z}[-2\htt(\nu)] )$:

\begin{multline} 
\gamma_{\{e\}}^{G_{\cO_n}}(\k_{T_\nu\cap Z}[-2\htt(\nu)] ) = Ra_!p^!i_*\k_{T_\nu\cap
Z}[-2\htt(\nu)] =
\\
R{\tii{a}}_!r^!\k_{T_\nu\cap Z}[-2\htt(\nu)] = R{\tii{a}}_!\k_{G_{\cO_n}\tim
(T_\nu\cap Z)}[2\dim (G_{\cO_n})-2\htt(\nu)] \,.
\end{multline} 
Let us consider a point $L_\eta\in Z$, where we again choose $\eta\in X_*(T)$ dominant. Then

\begin{lem} The dimension of the fiber $\tii{a}^{-1}(L_\eta)$ is $\dim G_{\cO_n}
-\htt(\nu+\eta)$ if $\eta\geq\nu$, otherwise it is empty. 
\end{lem}

\begin{proof}
We have
\begin{multline}
\tii{a}^{-1}(L_\eta) \ = \ \{(g,z)\in G_{\cO_n}\times (T_\nu\cap Z) \mid g\cdot z = \eta \} \ =
\\
\{(g,z)\in G_{\cO_n}\times (T_\nu\cap \GrG\eta) \mid g\cdot z = \eta \}\,.
\end{multline}
By theorem \ref{dimension estimates}, we see that $\tii{a}^{-1}(L_\eta)$ is empty unless
$\eta\geq\nu$. Furthermore, from theorem \ref{dimension estimates}, we see that if
$\eta\geq\nu$ then

\begin{multline}
\dim \tii{a}^{-1}(L_\eta) \ = \ 
\\
\htt(\eta-\nu) + \dim ((G_{\cO_n})_\eta) \ = \ 
\htt(\eta-\nu) +
\dim (G_{\cO_n}) -\dim \GrG\eta \ = \ 
\\
\htt(\eta-\nu) +
\dim (G_{\cO_n}) - 2\dim \htt (\eta)  \ = \ \dim (G_{\cO_n}) -\htt(\eta+\nu) \,.
\end{multline}

\end{proof}

In other words,
$\PPP_Z(\nu,\k)$ is the zero$^{th}$ perverse cohomology of the
$!$-image of a (shift of) a constant sheaf under an ``essentially semi-small'' map.
Here ``essentially semi-small'' means that the the dimensions of fibers 
have the correct increment but the generic fiber is not finite.

\section{\bf The structure of  projectives that represent weight functors}
\lab{The structure of  projectives that represent weight functors}

In this section we analyze the projective $ \PPP_Z(\k)=\pl_\nu\ \PPP_Z(\nu,\k)$ 
which represents the fiber functor $F$ on $\oP_\GO(Z,\k)$. As in the previous section, $Z$ is closed subset of $\G$ which is a union of finitely many $\GO$-orbits.

\begin{prop}
\label{structureProjectives}
(a) Let $Y\subset Z$ be a closed subset consisting of $\GO$-orbits. Then 
$$
\PPP_Y(\k)=\
{}^pH^0( \PPP_Z(\k)\rest Y)\,,$$
and there is a canonical surjection
$$
 \PPP_Z(\k) \to \PPP_Y(\k)\,.
$$

(b)
The projective $\PPP_Z(\k)$ has a filtration such that the the associated graded
$$
Gr(\PPP_Z(\k)) \ = \ \bb{\GG^\la\sub Z}\pl \ F[\II_*(\la,\k)]^*\ten_\k\II_!(\la,\k) \,.
$$
In particular, $F(\PPP_Z(\k))$ is free over $\k$. 

(c)
$ \PPP_Z(\nu,\k)\cong\
 \PPP_Z(\nu,\Z)
\bb{\Z}{\aa{L}\ten}
\k
$.

\end{prop}
\begin{proof} 
We begin with (a).  We write $i:\ Y\inj\ Z$ for the inclusion.  The identity
$\Hom(P_Z(\k),i_*-) =\Hom(i^*P_Z(\k),-)$ shows that the complex
$\PPP_Z(\k)\rest Y\in {}^pD^{\le 0}_\GO(Y,\k)$,
represents $F$ on the subcategory 
$\oP_\GO(Y,\k)$, and hence so does ${}^pH^0( \PPP_Z(\k)\rest Y)\in\oP_\GO(Y,\k)$. Thus, 
$\PPP_Y(\k)=\ {}^pH^0( \PPP_Z(\k)\rest Y)$. For any $\AA\in\oP_\GO(Y,\k)$ we have  the
identity $\Hom(P_Z(\k),\AA) =\Hom(P_Y(\k),\AA)$. This gives a canonical surjection
$\PPP_Z(\k)\to\PPP_Y(\k)$.

Now we will prove parts (b) and (c) simultaneously. We will argue (b) by induction on the
number of $\GO$-orbits in $Z$. Let us assume that $\GG^\la$ is an open $\GO$-orbit in
$Z$ and let $Y=Z-\GG^\la$. Let us consider the following exact sequence:
\begin{equation}\label{K}
0 \to K \to \PPP_Z(\k)\to \PPP_Y(\k)\to 0\,,
\end{equation}
where $K$ is simply the kernel of the canonical map from (a). Let $M$ be a $\k$-module
and let us  take
$\operatorname{RHom}$ of the exact sequence above to $\II_*(\la,M)$:

\begin{multline}
0\to
\Hom(\PPP_Y(\k),\II_*(\la,M)) \to \Hom(\PPP_Z(\k),\II_*(\la,M)) \to
\\
\Hom(K,\II_*(\la,M)) \to \Ext^1(\PPP_Y(\k),\II_*(\la,M))\,.
\end{multline}
By adjunction the first and last terms are zero and so we get, again by adjunction, 
\begin{multline}\label{ftor}
\Hom(K|_{\GG^\la},M_{\GG^\la}[2\htt(\la)])\cong\Hom(K,\II_*(\la,M)) \cong
\\
\Hom(\PPP_Z(\k),\II_*(\la,M)) = F(\II_*(\la,M))\,.
\end{multline}
We can view \eqref{ftor} as a functor from $\k$-modules to $\k$-modules
\begin{equation}
M \mapsto F(\II_*(\la,M))\,.
\end{equation}
This functor is, by the results in \S\ref{Canonical basis} represented by the free $\k$-module
$F(\II_*(\la,\k))^*$. As it is also represented by $K|_{\GG^\la}$, we conclude:

\begin{equation}
K|_{\GG^\la} \ = \ F(\II_*(\la,\k))^*\otimes_\k \k_{\GG^\la}[2\htt(\la)]\,.
\end{equation}

Now we claim:

\begin{equation}\label{KK}
K \ = \ F(\II_*(\la,\k))^*\otimes_\k \II_!(\la,\k) = \II_!(\la,F(\II_*(\la,\k))^*)\,.
\end{equation}

To prove this claim, let us consider the following exact sequence

\begin{equation}
0\to K' \to F(\II_*(\la,\k))^*\otimes_\k \II_!(\la,\k) \to K \to C \to 0\,.
\end{equation}
The kernel $K'$ and the cokernel $C$ are supported on $Y$. If we take
$\operatorname{RHom}$ of the exact sequence \eqref{K} to $C$, we get
\begin{equation}
0\to
\Hom(\PPP_Y(\k),C) \xrightarrow{\alpha} \Hom(\PPP_Z(\k),C) \to
\Hom(K,C) \to \Ext^1(\PPP_Y(\k),C)\,.
\end{equation}
Because $C$ is supported on $Y$ the map $\alpha$ is an isomorphism and the last term
vanishes. Therefore $C$ must be zero. To prove that $K'$ is zero, we first assume that
$\k=\Z$. Then, by \ref{!=!*}, $\II_!(\la,\Z)=\II_{!*}(\la,\Z)$. As $\II_{!*}(\la,\Z)$ has no
subobjects supported on $Y$, $K'=0$. Using \eqref{K} and \eqref{KK}, and proceeding by
induction on the number of $\GO$-orbits, we obtain (b) when $\k=\Z$. 

We will now prove (c). Because $\k\,\bb{\Z}{\aa{L}\ten}\II_!(\la,\Z)\cong\II_! (\la,\k)$
and because (b) holds for $\k=\Z$, we see that $\k\otimes_\Z \PPP_Z(\Z)$ is perverse. By
formula (2.6.7) in \cite{KS}, we see that for any $\AA\in \oP_\GO(Z,\k)$ we have
\begin{multline}
\operatorname{Hom}_\k(\k\otimes_\Z \PPP_Z(\Z), \AA) \ = \
\operatorname{Hom}_\Z(\PPP_Z(\Z), \operatorname{R}Hom(\k,\AA))\ = 
\\
\operatorname{Hom}_\Z(\PPP_Z(\Z), \AA) \ =  \ \mathbb {\oh}^*(\G,\AA)\,.
\end{multline}
Thus, $\k\otimes_\Z \PPP_Z(\Z)$ represents the functor $F$ on $\oP_\GO(Z,\k)$ and hence
we must have $\k\otimes_\Z \PPP_Z(\Z) =  \PPP_Z(\k)$. 

Finally to get statement (b) for an arbitrary ring $\k$, it suffices to use (c), (b) for the case
$\k=\Z$, and the fact that $\k\,\bb{\Z}{\aa{L}\ten}\II_!(\la,\Z)\cong\II_! (\la,\k)$.

\end{proof}

Let us write $\oP^{\k-proj}_\GO(Z,\k)$ for the subcategory of $\oP_\GO(Z,\k)$
consisting of sheaves $\cA\in\oP_\GO(Z,\k)$ such that ${\oh}^*(\GG,\cA)$ is
$\k$-projective. Note that, by lemma
\ref{perversity criterion}, the category
$\oP^{\k-proj}_\GO(Z,\k)$ is closed under Verdier duality. 
Because
$\PPP_Z(\k)\in\oP^{\k-proj}_\GO(Z,\k)$, its dual $I_Z(\k)=\D(\PPP_Z(\k))$ also
belongs in  $\oP^{\k-proj}_\GO(Z,\k)$. The sheaf
$I_Z(\k)$  is an injective object in the subcategory $\oP^{\k-proj}_\GO(Z,\k)$. Note
 that the abelianization of the exact category $\oP^{\k-proj}_\GO(\G,\k)$ is precisely $\pg$.

\section{\bf  Construction of the group scheme}
\label{Construction of the group scheme}

In this section we construct a  group scheme $\tiii G_\k$ such that $\pg$ is the category of its
representations. We proceed by Tannakian formalism, see, for example, \cite{DM}. Unlike
\cite{DM} we work over an arbitrary commutative ring $\k$. This is made possible by the
fact that $F(\PPP_Z(\k))$ is free over $\k$, \ref{structureProjectives}.

\begin{prop} \label{construction of tilde G}
There is a group scheme $\tiii G_\k$ over $\k$ such that the tensor category of
representations, finitely generated over $\k$, is equivalent to $\pg$.
Furthermore, the coordinate ring $\k[\tii G_\k]$ is free over $\k$ and
$\tiii G_\k = \operatorname{Spec}(\k)\times_ {\operatorname{Spec}(\Z)}\tiii
G_\Z$\,.
\end{prop}

\nop We view $\pg$ as a direct limit $\limi_Z\oP_\GO(Z,\k)$; here $Z$ runs through finite
dimensional $\GO$-invariant closed subsets of the affine Grassmannian $\G$. Let us
write $A_Z(\k)$ for the $\k$-algebra $\operatorname{End}(P_Z(\k)) = F(P_Z(\k))$.
The algebra $A_Z(\k)$ is free of finite rank over $\k$. Let us write
$\operatorname{Mod}_{A_Z(\k)}$ for the category of  $A_Z(\k)$-modules
which are finitely generated over $\k$. Because
$P_Z(\k)$ is a projective generator of $\oP_\GO(Z,\k)$, we see that the restriction of the
functor $F$ to $\oP_\GO(Z,\k) \ra \V$
lifts to an equivalence of abelian categories:
\begin{equation} 
\text{the categories $\oP_\GO(Z,\k)$ and $\operatorname{Mod}_{A_Z(\k)}$ are
equivalent as abelian categories.}
\end{equation}

As $A_Z(\k)$ is free of finite rank over $\k$, its $\k$-dual $B_Z(\k)$ is naturally a co-algebra.
Furthermore, let us consider a $\k$-module $V$. 
Because
\begin{equation}
\Hom_\k(A_Z(\k)\otimes_\k V, V) \ = \ \Hom_\k(V,B_Z(\k)\otimes_\k V)\,,
\end{equation}
we see that it is equivalent to give to $V$ a structure of an $A_Z(\k)$-module or to give it a structure of a
$B_Z(\k)$-comodule. Let us write $\operatorname{Comod}_{B_Z(\k)}$ for the category of  $B_Z(\k)$-comodules
which are finitely generated over $\k$.

From the previous discussion we conclude:
\begin{equation} 
\begin{gathered}
\text{the categories $\oP_\GO(Z,\k)$ and $\operatorname{Comod}_{B_Z(\k)}$}
\\
\text{are
equivalent as abelian categories.}
\end{gathered}
\end{equation}

Let us write $I_Z(\k)$ for the Verdier dual of $P_Z(\k)$. Now,
\begin{equation} 
B_Z(\k) \ = \ A_Z(\k)^* \ = \ \oh^*(\G,P_Z(\k))^* \ = \ \oh^*(\G,I_Z(\k))\ = \ F(I_Z(\k))\,.
\end{equation}
 
If  $Z \subset Z'$ are both closed and $\GO$-invariant then the canonical morphism
$P_{Z'}(\k) \to P_Z(\k)$ gives rise to a morphism $I_Z(\k) \to I_{Z'}(\k)$ and this, in turn,
gives a map of co-algebras $B_Z(\k) \to B_{Z'}(\k)$. Hence we can form the coalgebra
\begin{equation} 
B(\k) \ = \ \varinjlim B_Z(\k)\,,
\end{equation}
and we get:
\begin{equation} 
\text{the categories 
$\oP_\GO(\G,\k)$ and $\operatorname{Comod}_{B(\k)}$ are
equivalent as abelian categories.}
\end{equation}

It now remains to give the coalgebra $B(\k)$ the structure of an algebra and to give an inverse
in its coalgebra structure. We will start by giving $B(\k)$ an algebra structure. To
this end, let us consider the filtration of $\oP_\GO(\GG,\k)$ by the
subcategories
$\oP_\GO(\la,\k) =  \oP_\GO(\barr{\GG_\la},\k)$  indexed  dominant coweights $\la$. 
In the discussion that is to follow we use the following convention. When we substitute $\barr{\GG_\la}$ for the subvariety $Z$ we use the following shorthand notation
$A_\la(\k)=A_\Glb(\k),\ P_\la(\k)=P_\Glb(\k)$ etc.
This
filtration  is  compatible with the convolution product in the sense that $\oP_\GO(\la,\k)\ast
\oP_\GO(\mu,\k)\sub \oP_\GO(\la+\mu,\k)$. We have:
\begin{multline}\label{multiplication of projectives}
\Hom[P_{\la+\mu}(\k),P_{\la}(\k) \ast P_{\mu}(\k)]
\cong
F[P_{\la}(\k)
\ast P_{\mu}(\k)]
\cong
\\
F[P_{\la}(\k)]\tenk
F[P_{\mu}(\k)]
=
A_\lambda(\k) \otimes_\k A_\mu(\k)\,.
\end{multline}
The element 
$1\otimes 1 \in A_\la(\k) \otimes_\k A_\mu(\k)$ 
gives rise to a
morphism $P_{\la+\mu}(\k)\to P_{\la}(\k) \ast P_{\mu}(\k)$. Dualizing this gives a
morphism $I_{\la}(\k) \ast I_{\mu}(\k)\to I_{\la+\mu}(\k)$ and by applying the functor $F$ a
morphism
\begin{equation}
B_{\la}(\k) \otimes_\k B_{\mu}(\k)\to B_{\la+\mu}(\k)\,.
\end{equation}
Passing to the limit gives $B(\k)$ a structure of a commutative $\k$-algebra; the associativity
and the commutativity of the multiplication come from the associativity and commutativity of
the tensor product. To summarize, we have constructed an affine monoid $\tiii
G_\k = \Spec(B(\k))$ such that 
\begin{equation}
\operatorname{Rep}_{\tiii G_\k} \ \ \ \text{is equivalent to}\ \ \  \pg \ \ \ \text{as
tensor categories}\,;
\end{equation}
here $\operatorname{Rep}_{\tiii G_\k}$ denotes the category of representations
of $\tiii G_\k$ which are finitely generated as $\k$-modules. 

We will show next that $\tiii G_\k$ is a group scheme, i.e., that it has inverses. To do so, we first observe that
while the ind-scheme is $\GK$ is not of
ind-finite type it is a torsor
for the the pro-algebraic group
$\GO$, 
over two ind-finite type schemes
$\GK/\GO= \GG$ and $\GO\bss\GK$. This gives notions of
two kinds of equivariant perverse sheaves
on $\GG$
that come with equivalences
$
\oP_{\GO\tim 1}(\GK,\k)
\cong
\oP(\GO\bss\GK,\k)
$
and 
$\oP_{1\tim\GO}(\GK,\k)
\cong
\oP(\GG,\k)
$.
In particular one obtains two notions of full subcategories
$\oP_{\GO\tim\GO}(\GK,\k)$ of
$
\oP_{\GO\tim 1}(\GK,\k)
$
and 
$\oP_{1\tim\GO}(\GK,\k)
$, but these are easily seen to coincide.
Let us recall that we write $\oP^{\k-proj}_\GO(\G,\k)$ for the subcategory of $\oP_\GO(\G,\k)$
consisting of sheaves $\cA\in\oP_\GO(\G,\k)$ such that ${\oh}^*(\GG,\cA)$ is
$\k$-projective. The inversion map  $i:\GK\to\GK$, $i(g)=g^{-1}$
exchanges
$
\oP_{\GO\tim 1}(\GK,\k)
$
and 
$\oP_{1\tim\GO}(\GK,\k)
$
and defines an autoequivalence of
$\oP_{\GO\tim\GO}(\GK,\k)$
which we can view as
$i^*: \oP_\GO(\GG,\k)\to\oP_\GO(\GG,\k)$. 
Now we can define
an
anti-involution 
on $\oP^{\k-proj}_\GO(\G,\k)$  by
\begin{equation} 
\AA \mapsto \AA^*= \D(i^*\AA)\,.
\end{equation}
This involution makes $\operatorname{Rep}^{\k-proj}_{\tiii G_\k}$, the category
of representations of $\tiii G_\k$ on finitely generated projective $\k$-modules, a
rigid tensor category. By \cite{Sa} II.3.1.1, I.5.2.2, and I.5.2.3, we conclude that 
$\tiii G_\k = \Spec(B(\k))$ is a group scheme. 

 The  statement $\tiii G_\k =
\operatorname{Spec}(\k)\times_ {\operatorname{Spec}(\Z)}\tiii G_\Z$, i.e., that
$B(\k) = \k\otimes_\Z B(\Z)$,  now follows from Proposition
\ref{structureProjectives} part (c). The algebra $B(\k)$ is free over $\k$ by
construction.

\section{\bf The identification of $\tiii G_\k$ with the dual group of $G$}
\label{identification}

In this section we identify the group scheme $\tiii G_\k$. As $\tiii G_\k = 
\operatorname{Spec}(\k)\times_ {\operatorname{Spec}(\Z)}\tiii G_\Z$, by \ref{construction of
tilde G}, it suffices to do so when $\k=\bZ$. Recall that there exists a unique split reductive 
group scheme, the Chevalley group scheme, over $\bZ$ associated to any root datum (\cite{SGA},\cite{Dem}). 
Let $\cGZ,\cTZ$
be such schemes
associated to the root data dual to that of
$G,T$  and denote 
$\cGk=\cGZ\ten_\Z\k,\cTk=\cTZ\ten_\Z\k$ for any $\k$. 
We claim:

\begin{thm} \label{main result}
The group scheme $\tiii G_\bZ$ is the split reductive group scheme over $\bZ$  whose root
datum is dual to that of $G$. 
\end{thm}

The rest of this section is devoted to the proof of this theorem. We first recall that we have shown in section \ref{The case of a field of characteristic zero} that the above statement holds at the generic point, i.e., that $\tiii G_\bQ$ is split
reductive reductive group whose root datum is dual to that of $G$. By the uniqueness of the Chevalley group scheme, \cite{Dem}, and the fact that  $\tilde G_\bQ$ is a split
reductive reductive group whose root datum is dual to that of $G$ it suffices to show:

\begin{subequations} 
\begin{equation}\label{flat and smooth}
\text{The group scheme $\tiii G_\Z$ is  smooth over
$\operatorname{Spec}(\Z)$} 
\end{equation}
\begin{equation}\label{reductive at geometric points}
\text{At each geometric point $\operatorname{Spec}(\kappa)$, $\kappa=\bar
\bF_p$, the group scheme $\tiii G_\kappa$ is reductive}\,.
\end{equation}
\begin{equation} \label{maximal torus}
 \text{
The dual (split) torus $\check T_\bZ$ is a maximal torus of $\tiii G_\Z$}\,.
\end{equation}
\end{subequations} 

The properties \eqref{flat and smooth} and \eqref{reductive at geometric points} together with the fact that $\tiii G_\Z$ is affine amount to the definition of a reductive group. The last statement \eqref{maximal torus} says that $\tiii G_\Z$ is split. Note that it is not necessary to check in  \eqref{reductive at geometric points} that $\tiii
G_\kappa$ has root datum dual to that of $G$; that is automatic because it holds for  $\tiii G_\bQ$. However, to prove the fact that $\tiii G_\Z$ is smooth over $\operatorname{Spec}(\Z)$ and to deal with the fact that we do not yet know that $\tiii G_\Z$ is of finite type we end up having to calculate the root data of the $\tiii G_\kappa$. 

In what follows, we will make crucial use of results in \cite{PY}. We will first recall their theorem 1.5:

\begin{equation}\label{PY theorem 1.5}
\begin{gathered}
\text{An affine flat group scheme over the integers with all its fibers connected}
\\
\text{reductive algebraic groups of the same dimension is a reductive group\,.}
\end{gathered}
\end{equation}
As $\bZ[\tii G_\Z]$ is free over $\bZ$, by proposition \ref{construction of tilde G}, we see that $\tii G_\Z$ is flat. By section \ref{The case of a field of characteristic zero} all of the groups $\tii G_{\bQ_p}$ are split connected reductive groups with their root datum dual to that of $G$ and hence the $\tii G_{\bQ_p}$ are in particular connected reductive of the same dimension. Thus, to prove that $\tiii G_\Z$ is reductive, it suffices to show:
\begin{equation} \label{reductive over Z_p}
\text{The groups $\tii G_{\bZ_p}$ are reductive\,.}
\end{equation}
In order to prove \eqref{reductive over Z_p}, we will make use of the maximal torus, so we will now make a start at proving  \eqref{maximal torus}. We begin by exhibiting $\check T_\bZ$ as a sub torus of $\tiii G_\Z$. We proceed as in section 
\ref{The case of a field of characteristic zero}.
The  functor
\begin{equation} 
F ={\oh}^*: \operatorname P_{\GO}(\G,\bZ) \to \operatorname{Mod}_\Z(X_*(T)))
\end{equation} 
gives us a homomorphism 
$\cTZ \to \tGZ$. This makes  $\check T_\bZ$ a sub torus of $\tiii G_\Z$
because  for any cocharacter $\nu\in X_*(T)$, the $\nu$-weight space
$F_\nu(\II_!(\la,\Z))=H^*_c(S_\nu,\II_!(\la,\Z))$ is non-zero.

Let us now write $\kappa = \barr\bF_p$, for $p$ a prime, and $(\tGka)_{red}$ for the reduced subscheme of $\tGka$. We note that, just as in \ref{The case of a field of characteristic zero} we see that the group scheme $\tiii
G_\kappa$ is connected because  $\tii G_\kappa$ has no finite quotients -- there is no non-trivial
tensor subcategory of $\oP_\GO(\G, \barr\bF_p)$ supported on finitely many $\GO$-orbits. To complete the proof of Theorem  \ref{main result} we thus must argue, in addition to \eqref{reductive over Z_p}, that:
\begin{equation}
\text{The torus  $\check T_\kappa$ is maximal in $\tiii G_\kappa$\,.}
\end{equation}
We will argue these two points simultaneously. A crucial ingredient will be theorem 1.2 of \cite{PY}, which we now state in a form useful to us. The formulation below has much stronger hypotheses than in \cite{PY}. We inserted these stronger hypotheses as they are hold in case at hand to simplify the formulation. By theorem 1.2 of \cite{PY}, the group $\tii G_{\bZ_p}$ is reductive if the following conditions hold:

\begin{subequations}
\label{theorem 1.2 of PY}
\begin{equation} \label{assumption 1}
\text{\text{$\tii G_{\bZ_p}$  is affine and flat}}
\end{equation}
\begin{equation}  \label{assumption 2}
\text{$\tii G_{\bQ_p}$  is connected and reductive}
\end{equation}
\begin{equation} \label{assumption 3}
\text{$(\tGka)_{red}$ is a connected reductive group of the same type as  $\tii G_{\bQ_p}$\,.}
\end{equation}
\end{subequations}

As has been observed before, the first two hypotheses above are satisfied. Thus, it remains to prove \eqref{assumption 3}. To summarize, we are reduced to showing:
\begin{subequations}
\label{first reduction}
\begin{equation}
\label{first reduction a}
\text{the torus $\check T_\kappa$ is maximal in $(\tiii G_\kappa)_{red}$ and}
\end{equation}
\begin{equation}
\label{first reduction b}
\text{the group scheme
$(\tiii G_\kappa)_{red}$ is reductive with root datum
dual to that of $G$}\,.
\end{equation}
\end{subequations}
The two statements above are statements about the fiber at $\kappa=\barr \bF_p$. In the rest of this section we will be working at this fiber, but before doing so, we make one more
argument using the entire flat family. The flatness of
$\tilde G_{\bZ_p}$ implies that
\begin{equation} \label{flatness and fiber dimension}
\dim G = \dim \tiii G_{\bQ_p} \geq \dim (\tiii G_\kappa)_{red}\,.
\end{equation}
To see this, let us choose algebraically independent elements in the coordinate ring
of $(\tiii G_{\bF_p})_{red}$ and lift them to the coordinate ring of $\tiii G_{\bZ_p}$. By flatness,  $\bZ_p[\tiii G_{\bZ_p}] \subset\bQ_p[\tiii G_{\bQ_p}]$ and hence these lifts remain algebraically independent in the coordinate ring of $\tiii G_{\bQ_p}$. This gives  \eqref{flatness and fiber dimension}. Note that we do not a priori have an equality in \eqref{flatness and fiber dimension}, as we do not yet know that $\tilde G_{\Z_p}$ is of finite type.   

Now we are ready to work at $\kappa$. We proceed at the beginning in the same way as we did in section \ref{The case of a field of characteristic zero}. 
First, we have:
\begin{equation}
\label{irreducible representations in finite characteristic}
\begin{gathered}
\text{The irreducible representations of $\tiii G_\kappa$ are}
\\
\text{parameterized by dominant
weights $\la\in X^*(\check T_\kappa) = X_*(T)$\,.}
\end{gathered}
\end{equation}
We see this, just as in section \ref{The case of a field of characteristic zero}, by noting that the irreducible objects in $\operatorname P_{\GO}(\G,\barr\bF_p)$ are given by the $\II_{!*}(\la,\barr\bF_p)$, for $\la\in X_*(T)$ dominant.
Let us write, as in section \ref{The case of a field of characteristic zero}, $L(\la)$ for the irreducible representation of $\tiii G_\kappa$ associated to $\la$. First of all, because of proposition  \ref{CanonicalBasis}, we see that the weights of the representation $W(\lambda,\k)$ corresponding to $\II_{!}(\la,\k)$ are independent of $\k$. Hence, those weights are precisely the weights of the irreducible representation of $\check G_\bC$ of highest weight $\lambda$. On the other hand, we can write $\II_{!*}(\la,\barr\bF_p)$ in the Grothendieck group as a sum involving $\II_{!}(\la,\barr\bF_p)$  and terms $\II_{!}(\mu,\barr\bF_p)$, for $\mu < \la$. Hence, we conclude:
\begin{subequations}\label{structure of representations in positive characteristic}
\begin{equation}
\begin{gathered}
\text{The $\check T_\kappa$-weights of the irreducible
representation $L(\la)$ are contained}
\\
\text{in the  $\check T$-weights  of the irreducible representation of $\check G_\bC$}
\\
\text{associated to $\la$, and $\la$ is the highest weight in  $L(\la)$}\,.
\end{gathered}
\end{equation}
\begin{equation}
\text{The weights of $L(\lambda)$ are symmetric under the Weyl group $W$.}
\end{equation}
\begin{equation}\label{structure of weights in positive characteristic}
\begin{gathered}
\text{For $\alpha$ a simple positive root of $G$, with the corresponding coroot}
\\
\text{$\check\alpha$, the weights of $L(\lambda)$ on the line segment between}
\\
\text{$s_{\check \alpha}(\lambda)$ and $\lambda$ are all of the form $\lambda-k\check\alpha$, with $k$ an integer.}
\end{gathered}
\end{equation}
\end{subequations}
Note that, contrary to the case of characteristic zero, not all the weights between $\la$ and $\lambda-k\check\alpha$ occur as weights of $L(\la)$.

Next, we approximate
$\tGka$ by finite type quotients
$\sG$.
For any group scheme $H$ 
let us write $\operatorname{Irr}_H$ for the set of
irreducible representations of $H$.  We choose a  quotient group scheme
$\sG$ of 
$\tGka$ with the following properties:

\begin{subequations}
\label{inverse limit of group schemes of finite type over a field}
\begin{equation}
\text{$\sG$ is of finite type}
\end{equation}
\begin{equation} 
\label{equality of irreducibles}
\text{the canonical map 
$\operatorname{Irr}_{\sG} \to
\operatorname{Irr}_{\tGka}$ is a bijection}
\end{equation}
\begin{equation}
(\sG)_{red}=(\tGka)_{red}
\end{equation}
\end{subequations}

It is possible to satisfy the first and third conditions because any group scheme is a
projective limit of group schemes of finite type and by \eqref{flatness and fiber
dimension} the group scheme
$(\tGka)_{red}$ is of finite type. 
To ensure that 
\eqref{equality of irreducibles} is satisfied, 
it is enough
to choose 
$\sG$ sufficiently large so that the irreducible representations 
$L(\la)$ associated to a finite set of generators $\la$
of the semigroup of dominant cocharacters,
are pull-backs of representations of 
$\sG$. For any dominant
coweights $\la,\mu$, the sheaf 
$\II_{!*}(\la+\mu, \barr\bF_p)$ is  a subquotient
of the convolution product $\II_{!*}(\la, \barr\bF_p)\ast\II_{!*}(\mu, \barr\bF_p)$
since the support of the convolution is  
$\barr\GG_{\la+\mu}$. This shows that  all
irreducible representations 
$L(\la)$ of $\tGka$ come from
$\sG$.

Let us write $R$ for the reductive quotient of $(\tiii
G_\kappa)_{red}=(\sG)_{red}$ and note that, of course, $\check T_\kappa$ lies naturally in $R$. As any irreducible representation of
$(\sG)_{red}$ is trivial on the unipotent radical we have:
\begin{equation}
\text{The canonical map $\operatorname{Irr}_{R} \to\operatorname{Irr}_{(\sG)_{red}}
$ is a bijection\,.}
\end{equation} 

We now argue that in order to prove \eqref{first reduction} it suffices to
show that:
\begin{subequations}
\label{second reduction}
\begin{equation}
\label{second reduction a}
\text{the torus $\check T_\kappa$ is maximal in $R$}
\end{equation}
\begin{equation}
\label{second reduction b}
\text{the  root datum of $R$ with respect to $\check T_\kappa$ is dual to that of $G$}\,.
\end{equation}
\end{subequations}
Let us, then, assume \eqref{second reduction}.  We conclude immediately that  $\dim(R) = \dim(G)$, and, together with 
\eqref{flatness and fiber dimension}, this gives
\begin{equation}
\dim G = \dim \tilde G_\mathbb Q
\ge \dim (\sG)_{red} \geq \dim R =\dim G\,.
\end{equation}
Thus, we must have $(\sG)_{red}=R$ and hence $(\sG)_{red}$
is reductive  with its root datum dual to that of $G$. 
Since this holds for each of approximation
$\sG$ of
$\tGka$ we see that
$(\tGka)_{red}$ coincides with $R$, has
$\check T_\kappa$ as its  maximal torus,  and   its root datum dual to that of $G$. This gives 
\eqref{first reduction}.

The proof of
\eqref{second reduction} will be based on 
relating representations of $\tGka$
and $R$ (i.e., of $\sG$ and $(\sG)_{red}$),
by considering the $n\thh$  
powers of the
Frobenius maps between the $\ka$-scheme 
$\sG$ and  its $n\thh$ Frobenius 
twist $(\sG)^{(n)}$,  as depicted in the diagram below: 
\begin{equation}
\begin{CD}
\sG
@>{\operatorname{Fr}^n_{\sG}}>>
(\sG)^{(n)}
\\
@AAA @AAA
\\
 (\sG)_{red}
@>{\operatorname{Fr}^n_{\sG}}>>
(\sG)_{red}^{(n)}
.
\end{CD}
\end{equation} 
Since $\sG$ is of finite type we see, by
\cite[corollary III.3.6.4]{DG}, that
it is isomorphic to $(\sG/(\sG)_{red})\tim (\sG)_{red}$
as a scheme with the right 
multiplication action by 
$(\sG)_{red}$,\
and that the coordinate ring of
$\sG/(\sG)_{red}$ is of the form
$
\kappa[X_1,...,X_n]/\lb
X_1^{p_1},...,X_n^{p_n}\rb
$ for some powers
$p_i=p^{e_i}$ of $p$.
Hence, for $n\ge e_i$,  
\begin{equation}
\text{the Frobenius map
$\operatorname{Fr}^n: 
\sG
\to (\sG)^{(n)}$ factors through
$(\sG)_{red}^{(n)}$}\,.
\end{equation} This implies:
\begin{equation}
\begin{gathered}
\text{the $\text{n}^{th}$
 Frobenius twists of irreducible representations}
\\
\text{of 
$(\sG)_{red}$ extend to $\sG$
}\,.
\end{gathered}
\end{equation}

From the previous discussion we conclude, first of all, that
\begin{equation}
\text{We have an injection $\operatorname{Irr}_{R^{(n)}}\inj\ 
\operatorname{Irr}_{\tGka}$\,.}
\end{equation} 
This means that the torus $\cTka$ must be maximal in $R$, as the statement above bounds the size of the lattice of irreducible representations of $R$. Furthermore, we also conclude that 
\begin{equation}\label{equality of some irreducibles}
L^{\tGka}(p^n\la) \ = \ L^{R}(p^n\la) \ \ \ \text{for all $\la\in X_*(T)$ dominant}\,.
\end{equation}
We will now argue the second part of \eqref{second reduction}, i.e., that :
\begin{equation}\label{equality of root data in positive characteristic}
\text{
$
\De(R, \check T_\kappa) = \check \De(G,T)
$ 
\ \ \ and \ \ \ 
$
\check \De(R,\check T_\kappa) =\De(G,T)
$
}\,. 
\end{equation}
The argument here is a bit more involved than the argument in characteristic zero in
section \ref{The case of a field of characteristic zero}, but the basic idea is the same:
the pattern of the weights of irreducible representations determines the root datum.

The statement \eqref{structure of representations in positive characteristic} expresses the patterns of weights of the representations $L^{\tGka}(\la)$. The pattern of weights of the $L^{R}(\la)$ has a similar description, as $R$ is a reductive group. Comparing these patterns for $p^n\la$, we conclude that that the walls of the Weyl chambers of the root systems  of $(G,T)$ and $(R,\check T_\kappa)$ coincide in $X_*(T)$, and furthermore,  that the simple root directions of $R$ coincide with the simple coroot directions of $G$. Recall that in  characteristic zero we obtained an
equality of simple roots of $R$ and the simple coroots of $G$, but we cannot immediately conclude this fact here, as we only have a containment in \eqref{structure of weights in positive characteristic}. Hence, we must argue further.

As the next step,  we prove two
inclusions of lattices:
\begin{equation}\label{inclusions of lattices}
\Z\cd\De(R,\check T_\kappa) \ \subset \ \Z\cd\check\De(G,T)
\ \ 
\text{and}\ \
\Z\cd\De(G,T) \ \subset \ \Z\cd\check\De(R,\check T_\kappa)\,.
\end{equation}
We do so by analyzing the centers. First of all, note that 
the center of the reductive  group $R$ can be identified with the group scheme
\begin{equation}
\Hom(X^*(\ch T_\ka)/\Z\cd \De(R,\check
T_\kappa),\bG_{m,\kappa})\,.
\end{equation}
On the other hand, the tensor category 
$\PP_\GO(\GG,\barr\bF_p)$ is naturally equipped with a  grading by the
abelian group  
$\pi_0(\GG)\cong \pi_1(G)
=X_*(T)/\Z\cd\check\De(G,T)$, the group of connected components of $\GG$. Note that this grading is compatible with the tensor structure. This grading exhibits the
group  scheme
\begin{equation}
\Hom(X_*(T)/\Z\cd   \check\De(G,T), \bG_{m,\kappa})
\end{equation}
 as a group subscheme of the center
of  $\tGka$. Since this group subscheme  lies in $\cTka$
it also lies in the center of $R$.
This inclusion of centers corresponds to a natural surjection
\begin{equation}
X^*(\ch T_\ka)/\Z\cd \De(R,\check T_\kappa) \ \ 
\ \ \twoheadrightarrow \ \  X_*(T)/\Z\cd  
\check\De(G,T)\,.
\end{equation}
This gives the first inclusion in \eqref{inclusions of lattices}.
Since the lattices have bases of simple (co)roots
which have the same directions
the second inclusion follows from
$\left<\check\alpha,\alpha\right>=2$. Let us note that by the same reasoning the proof of \eqref{equality of root data in positive characteristic} can now be reduced to proving either of the equalities in \eqref{inclusions of lattices}. We will do so by building up from special cases.

First, let us observe that we are done when 
$G$ is of adjoint type. In that case
$\Z\cd\De(G,T)= X^*(T)$ and so we must have $\Z\cd\De(G,T)= 
\Z\cd\check\De(R,\cTka)$. At this point we no longer have need to argue in terms of $R$. In what follows we will prove directly that $\tii G_\kappa$ is reductive with its root datum dual to that of $G$. 

Next,  assume  that $G$ is
semi-simple and write
\begin{equation}\label{G to the adjoint group}
G \twoheadrightarrow G_\text{ad}\,, \ \ \ \text{where $G_\text{ad}$ is the adjoint quotient of $G$}\,.
\end{equation}
We have already shown that $(\tii G_\text{ad})_\kappa\cong (\check G_\text{ad})_\kappa$. Just as above, we note that we can provide the category $\oP_{G_{\text{ad},\cO}}(\G_{G_\text{ad}},\barr\bF_p)$ with a grading, as a tensor category, by the finite group $\pi_1(G_\text{ad})/\pi_1(G)$. With this grading the tensor subcategory 
$\oP_{G_{\cO}}(\G_{G},\barr\bF_p)$ of
$\oP_{G_{\text{ad},\cO}}(\G_{G_\text{ad}},\barr\bF_p)$ corresponds to the identity coset $\pi_1(G)$. Thus, we obtain  a 
surjective homomorphism
\begin{equation}\label{dual of G to the adjoint group}
(\check G_\text{ad})_\kappa = (\tii G_\text{ad})_\ka \twoheadrightarrow \tGka
\end{equation}
with a finite central kernel, given precisely by $\Hom(\pi_1(G_\text{ad})/\pi_1(G), \bG_{m,\kappa})$.
This implies that $\tGka$ is reductive. Let us write $T_\text{ad}$ for the maximal torus in  $G_\text{ad}$. Then from \eqref{G to the adjoint group} we see that the roots and the coroots of the pairs $(G,T)$ and $(G_\text{ad},T_\text{ad})$ coincide under the surjection $T\to T_\text{ad}$ and similarly, from \eqref{dual of G to the adjoint group}, we conclude that the roots and coroots of the pairs $(\tGka, \cTka)$ and $((\check G_\text{ad})_\kappa, (\check T_\text{ad})_\kappa)$ coincide under the surjection $ (\check T_\text{ad})_\kappa \to \cTka$. Thus, as we know the result for the adjoint group, we conclude that the root datum of $(\tGka, \cTka)$ is dual to that of $(G,T)$.

Finally, consider the case of a general reductive $G$. 
Let us write $S=Z(G)^0$ for the connected
component of the center of $G$. Then we have an exact sequence
\begin{equation}
1 \to S \to G \to G_\text{der}\to 1\,,
\end{equation}
where the derived group $G_\text{der}$ of $G$ is semisimple. 
This gives  maps:
\begin{equation}
\begin{CD}
\G_S @>i>> \G_G @>{\pi}>> \G_{G_\text{der}}\,,
\end{CD}
\end{equation}
which exhibit $\G_G$ as a trivial cover of $\G_{G_\text{der}}$ with fiber 
$\G_S$. By taking pushfowards
of sheaves this gives us the following sequence of functors:
\begin{equation}
\begin{CD}
\op_{S_\cO}(\G_S,\barr\bF_p) @>{\omega}>> \op_{G_\cO}(\G_G,\barr\bF_p)
@>{\gamma}>>\op_{{G_\text{der}}_\cO}(\G_{G_\text{der}},\barr\bF_p) \,, 
\end{CD}
\end{equation}
where
$\omega$ is clearly an embedding and $\gamma$ is essentially surjective because of the triviality
of the cover.  This, in turn, gives the following exact sequence of
group schemes:
\begin{equation}
1 \to 
(\tii G_\text{der})_\ka \to 
\tGka \to \tii S_\ka \to 1
\end{equation}
The fact that we have exactness at both ends follows from the fact $\omega$ is an embedding and
$\gamma$ is essentially surjective. To see the exactness in the middle, let us consider the quotient 
$\tGka/(\tii  G_\text{der})_\ka$. 
The representations of this groups scheme are given by objects in
$\op_{G_\cO}(\G_G,\ka)$ whose push-forward under $\pi$ to 
$\G_{G_\text{der}}$ consists of direct sums
of trivial representations. 
But these constitute precisely $\op_{S_\cO}(\G_S,\ka)$. Thus, $\tGka$ is
reductive. Arguing just as in the previous step, using the fact that the roots and coroots of $G$ and $G_\text{der}$ on one hand and those of $\tGka$ and $(\tii  G_\text{der})_\ka$ on the other,  coincide, we conclude that the root datum of $\tGka$ is dual to that of $G$.

\section{\bf Representations of reductive groups}
\lab{Representations of reductive groups}

The point of view we have taken in this paper so far is that of giving a canonical, geometric
construction of the dual group. In this section we turn things around and view our work as giving a
geometric interpretation of representation theory of split reductive groups. 

As before, $\k$ is commutative ring, noetherian and of finite global dimension. Recall the content of our 
main theorem \ref{main result}:
\begin{equation}\label{restatement}
 \operatorname{Rep}_{\check G_\k}\ \  \text{is equivalent to} \ \ \pg\ \ \text{as tensor categories};
\end{equation}
here $\operatorname{Rep}_{\check G_\k}$ stand for the category of $\k$-representations of $\check G_\k$,
finitely generated over $\k$, and $\check G_\k$ stands for the canonical split group scheme associated to the
root datum dual to that of the complex group $G$. This way we get a geometric interpretation of
representation theory of $\check G_\k$. The case when $\text{Char}(\k)=0$ was discussed in section
\ref{The case of a field of characteristic zero}.

Following our previous discussion we have $\check T_\k \subset \check B_\k \subset \check G_\k$, a maximal
torus and a Borel in $\check G_\k$. Associated to a weight $\la\in X_*(\check T)$ there are two standard
representations of highest weight $\lambda$. Let us describe these representations. We extend $\la$ to a
character on $\check B_\k$ so that it is trivial on the unipotent radical of  $\check B_\k$ and then induce this
character to a representation of $\check G_\k$. We call the resulting representation the Schur module and
denote it by $S(\lambda)$. As a module it is free over $\k$. The other representation associated to $\lambda$ is
the Weyl module $W(\lambda) = S(-w_0\lambda)^*$, where $w_0$ is the longest element in the Weyl group.
There is a canonical morphism  $W(\lambda) \to S(\lambda)$ which is the identity on the $\lambda$-weight
space. We have:
\begin{prop}
Under the equivalence \eqref{restatement} the diagrams
$W(\lambda) \to S(\lambda)$ and  $\mathcal I_!(\la) \to \mathcal I_*(\la)$ correspond to each other. 
\end{prop}

\begin{proof}
The modules  $S(\lambda)$ and $W(\lambda)$ can also be characterized in the following manner. Let us
write $\operatorname{Rep}_{\check G_\k}^{\leq \la}$ for the full subcategory of
$\operatorname{Rep}_{\check G_\k}$ consisting of representations whose $\check T_\k$-weights are all
$\leq\la$. Then the representations
$S(\la)$ and $W(\la)$ satisfy the following universal properties:
\begin{subequations}\label{universal properties}
\begin{equation}
\text{for $V\in\operatorname{Rep}_{\check G_\k}^{\leq \la}$ \ \ we have \ \ $\Hom_{\check
G_\k}(V,S(\la))=\Hom_\k(V_\lambda,\k)$}
\end{equation}
and
\begin{equation}
\text{for $V\in\operatorname{Rep}_{\check G_\k}^{\leq \la}$ \ \ we have \ \ $\Hom_{\check
G_\k}(W(\la),V)=\Hom_\k(\k, V_\la)$}
\end{equation}
\end{subequations}
On the geometric side the category $\operatorname{Rep}_{\check G_\k}^{\leq \la}$ corresponds to the
category $\oP_\GO(\la,\k) =  \oP_\GO(\barr{\GG_\la},\k)$. Obviously, the sheaves $\mathcal I_*(\la)$ and
$\mathcal I_!(\la)$ belong $\oP_\GO(\la,\k)$ and satisfy the universal properties \eqref{universal properties},
proving the proposition.
\end{proof}

As a corollary, proposition \ref{CanonicalBasis} gives:

\begin{cor}
The $\nu$-weight spaces $S(\la)_\nu$ and $W(\la)_\nu$ of $S(\la)$ and $W(\la)$, respectively, can both be
canonically identified with the $\k$-vector space spanned by the irreducible components of
${\barr{\GG_\la}}\cap S_\nu$.
In particular, the dimensions of these weight spaces are given 
by the number of irreducible components of ${\barr{\GG_\la}}\cap S_\nu$\,.
\end{cor}

Finally, let us assume that $\k$ is a field. Then we also have an irreducible representation $L(\la)$ associated
to $\la$. Under the correspondence  \eqref{restatement} the irreducible representation $L(\la)$ corresponds to
the irreducible sheaf $\mathcal I_{!*}(\la,\k)$.

An important motivation for our work was its potential application to representation theory of algebraic groups. To this end, we would like to propose the following:
\begin{conj}
The stalks of $\mathcal I_!(\la,\bZ)$ are free.
\end{conj}

\section{\bf Variants and the geometric Langlands program}
\lab{Variants and the geometric Langlands program}

As was stated before, in this paper we have worked with $\bC$-schemes because in this case we have a good
sheaf theory for sheaves with coefficients in any commutative ring, in particular, the integers. It is also possible
to work with other topologies. This is important for certain applications, for example for the geometric
Langlands program, since our results can be viewed as providing the unramified local geometric Langlands
correspondence. 

We will explain briefly the modifications necessary to work in the etale topology and over an arbitrary
algebraically closed base field $K$. To this end, let us view the group $G$ as a
split reductive group over the integers. All the geometric constructions made in this paper go through over the
integers, in particular, our Grassmannian
$\G$ is defined over $\bZ$. We write $\G_K$ for the affine Grassmannian over the base field $K$. In a few
places in the paper we have argued using
$\bZ$ as coefficients, for  instance, in section \ref{identification}. When we work in the etale topology, we simply
replace $\bZ$ by $\bZ_\ell$, where $\ell\neq Char(K)$. 

For completeness, we state here a version of our theorem for $\G_K$:

\begin{thm} 
There is an equivalence of tensor categories
$$
\op_{G(\Cal O)}(\G_{K},\k	) \cong Rep(
\check G_{	\k	}	).
$$
where we can take $\k$ to be any ring for which the left hand side is defined and which can be obtained by
base change from $\bZ_\ell$, for example, $\k$ could be $\bQ_\ell$, $\bZ_\ell$, $\bZ/\ell^n\bZ$, $\barr {\bF_\ell}$. 
\end{thm} 

\begin{rmk}
The previous theorem allows one to extend the notion of Hecke
eigensheaves  in the geometric Langlands program from the case of characteristic zero coefficients to coefficients in an arbitrary field. This is used in \cite{G} which gives a proof of de Jong's conjecture. 
\end{rmk}

\appendix

\section{\bf
Categories of  perverse sheaves}
\label{Categoryofequivariantperversesheaves}

In this appendix we prove propositions \ref{equivalence} and
\ref{another equivalence}, i.e., we will show that

\begin{prop}
The categories $\op_{\GO \rtimes \operatorname{Aut}(\OO)} (\G,\k)$, 
$\pg$, and $\ps (\G,\k)$ are naturally
equivalent.
\end{prop}

We first note that  $\op_{\GO \rtimes \operatorname{Aut}(\OO)} (\G,\k)$ is a full
subcategory of $\pg$ and $\pg$ is a full subcategory of $\ps (\G,\k)$; this follows
from the fact that stabilizers of points are connected for the actions of $\GO \rtimes
\operatorname{Aut}(\OO)$ and $\GO$ on $\G$. The reductive quotients of $\GO
\rtimes \operatorname{Aut}(\OO)$, $\GO$, and $\operatorname{Aut}(\OO)$ are
$G\times \bG_m$, $G$, and $\bG_m$, respectively. Hence, it suffices to show
that for any $\GO$-invariant finite dimensional subvariety we have:
\begin{equation} \label{equivalences}
\text{the categories $\ \op_{\cS,G} (Z,\k)\ $ and $\ \op_{\cS,\bG_m} (Z,\k)\ $ are
equivalent to $\ \ps (Z,\k)$}\,;
\end{equation}
here $\op_{\cS,G} (Z,\k)$ and $\op_{\cS,\bG_m} (Z,\k)$ denote subcategories
of $\ps (Z,\k)$ consisting of sheaves which are $G$ and $\bG_m$- equivariant,
respectively. 

We will proceed by induction in the following manner. Obviously, the
statement above holds if $Z$ is a $\GO$-orbit. Hence, by induction, it suffices to
prove  \eqref{equivalences} for a $\GO$-invariant subset $Z$ under the following
hypotheses: 
\begin{equation}\label{induction step}
\begin{gathered}
\text{for some dominant $\lambda$, the orbit $\G^\lambda$ is closed in 
$Z$}
\\
\text{ and \eqref{equivalences} holds for the open set $U=Z-\G^\lambda$}\,.  
\end{gathered}
\end{equation}

To prove the above statement, we use the gluing construction of \cite{MV},
\cite{V} for perverse sheaves. Let us recall the construction. We write 
$\cA$ and $\cB$ be two abelian categories and
$F_1,G_1:\cA
\to
\cB$ two functors, $F_1$ right exact, $F_2$ left exact, and $T:F_1\to F_2$ a
natural transformation. We define a category $\cC(F_1,F_2;T)$ as follows. The
objects of $\cC(F_1,F_2;T)$ consist of pairs of objects $(A,B) \in
\text{Ob}(\cA) \times \text{Ob}(\cB)$ together with a factorization
$F_1(A) \xrightarrow{m} B \xrightarrow{n} F_2(A)$ of $T(A)$, i.e., $n\circ m
= T(A)$. The morphisms of  $\cC(F_1,F_2;T)$ are given by pairs of morphisms
$(f,g) \in
\text{Mor}(\cA) \times \text{Mor}(\cB)$ which make the appropriate
prism commute. The category $\cC(F_1,F_2;T)$ is abelian.

We use this formalism in various situations. To begin with, let us write
$j:U\hookrightarrow Z$ for the inclusion and set:
\begin{equation}
\begin{split}
&\cA = {\ps}(U,\k)
\\
&\cB =\V
\\
&F_1 = F_\la \circ {^pj_!}
\\
&F_2 = F_\la \circ {^pj_*}
\\
&T= F_\la( {^pj_!} \to {^pj_*})\,.
\end{split}
\end{equation}

We have a functor
\begin{equation}
E: \ps(Z,\k) \to \cC(F_1,F_2;T)
\end{equation}
which sends $\cF\in\ps(Z,\k)$ to $A = \cF \rest  U$, $B = F_\la(\cF)$ and
the factorization $F_1(A) \xrightarrow{m} B \xrightarrow{n} F_2(A)$ is the
one gotten by applying
$F_\la$ to ${^pj_!}(\cF\rest U) \to \cF \to {^pj_*}(\cF\rest U)$. By Proposition 1.2
in \cite{V} the functor $E$ is an embedding. Two remarks are in order. First, in
\cite{V} we work over a field, but this is not used in the proof. Secondly, the
functor $E$ is actually an equivalence of categories. 

Let us now bring in the group $G$. We write $\tii\cS$ for the stratification of
$G\times Z$ by subvarieties $G\times \G^\la$. We write $a:G\times Z \to Z$ for the
action map and $p:G\times Z\to Z$ for the projection. Let $\cF\in \ps(Z,\k)$. We
have an isomorphism $\phi:p^*\cF|G\times U\cong a^*\cF|G\times U$ such that
$\phi|\{e\}\times U = \text{id}$. We are now to extend the $\phi$ to $G\times Z$.
To this end we first construct a functor $\tii 
F_\la: \op_{\tii\cS}(G\times Z,\k) \to \tii\cB$, where $\tii\cB$ stands for the
category of $\k$-local systems on $G$. Let us write $\tii S_\la = G\times
S_\la$, denote by $i:\tii S_\la \hookrightarrow G\times Z$ the inclusion, and
write $\pi:G\times Z \to G$ for the projection. Then $\tii F_\la = R\pi_!i^*$.
Furthermore, we write $j:G\times U \hookrightarrow G\times G$ for the inclusion
and set:

\begin{equation}
\begin{split}
&\tii\cA = {\ps}(G\times U,\k)
\\
&\tii \cB =\{\text{$\k$-local systems on $G$}\}
\\
&\tii F_1 = \tii F_\la \circ {^p\tii j_!}
\\
&\tii F_2 = \tii F_\la \circ {^p\tii j_*}
\\
&\tii T= \tii F_\la( {^p\tii j_!} \to {^p\tii j_*})\,.
\end{split}
\end{equation}
As before, we have a functor
\begin{equation}
\tii E: \op_{\tii \cS}(G\times Z,\k) \to \cC(\tii F_1,\tii F_2;\tii T)
\end{equation}
which sends $\tii\cF\in\ps(G\times Z,\k)$ to $\tii A = \tii\cF \rest 
G\times U$,
$\tii B = \tii F_\la(\tii\cF)$ and the factorization $\tii F_1(\tii A)
\xrightarrow{m} \tii B
\xrightarrow{n} \tii F_2(\tii A)$ is the one gotten by applying
$\tii F_\la$ to ${^p\tii j_!}(\tii\cF\rest G\times U) \to \tii\cF \to
{^p\tii j_*}(\tii\cF\rest G\times U)$. By the same reasoning as above, the
functor $\tii E$ is an equivalence of categories. 

We will now apply $\tii E$ to $p^*\cF$ and to $a^*\cF$. For $\tii E (p^*\cF)$
we get the data of $E(\cF)$ at $\{e\}\times Z$ extended across $G$ as the constant
local system. For $\tii E (a^*\cF)$ we also get the extension data of $E(\cF)$ at
$\{e\}\times Z$. Because, by theorem \ref{exactness}, the functors $F_\nu$ are
independent of the data $T\subset B$ used in defining them, we see that the
extension data for $a^*\cF$ restricted to $\{g\}\times Z$, for any $g\in G$, is
canonically identified with the extension data of $a^*\cF$ at $\{e\}\times Z$. This
gives us an identification of $\tii E (a^*\cF)$ with $\tii E (p^*\cF)$ and hence
an isomorphism between $a^*\cF$ and $p^*\cF$. This shows the first part of
\ref{equivalences}. 

The case of the group $\bG_m$ is even a bit simpler. Here we use the fact that
$\bG_m$-action preserves the variety $S_\la$ and hence all the
$\bG_m$-translates of the functor $F_\la$ are identical to $F_\la$. 

\section{\bf Corrections to section 12}

In section 12 of the paper  we approximate $\tGka$ by finite type quotients $\sG$ satisfying properties (12.11). One should not have  imposed condition (12.11c) as at this point we do not yet know if any finite type quotients $\sG$ satisfying it exist. In the text it is erroneously  claimed that (12.8) implies that $(\tGka)_{red}$ is of finite type when in reality (12.8) only asserts that $(\tGka)_{red}$ is of finite dimension. However, the argument in the paper goes through without imposing (12.11c) as follows. We write $\tGka = \varprojlim \sG$ with the $\sG$ satisfying (12.11a) and (12.11b). Then we have $(\tGka)_{red} = \varprojlim (\sG)_{red}$. At the end of section 12 it is shown that all the  $(\sG)_{red}$ are connected reductive and isomorphic to $\check G_\kappa$. If we have two finite type quotients $\sG$ and $\tGka^{**}$ of $\tGka$ satisfying (12.11a) and (12.11b) such that $\tGka \to \tGka^{**}\to \sG$ then the induced map $ (\tGka^{**})_{red} \to (\sG)_{red}$ amounts to an endomorphism of $\check G_\kappa$ which is an identity on the maximal torus $\check T_\kappa$ and a bijection on the irreducible representations preserving the highest weights. Such an endomorphism is a unipotent isogeny and, for example, by~\cite[2.3 corollary]{PY} is an isomorphism. Thus we conclude that $(\tGka)_{red} = \check G_\kappa$ which proves (12.6) and thus completes the identification of the group scheme $\tilde G_\bZ$ with $\check G_\bZ$.

\end{document}